\documentclass[12pt]{article}

\usepackage{graphicx}
\usepackage{latexsym,amsmath,amsfonts,amscd, amsthm}
\usepackage{epsfig}
\usepackage{changebar}
\usepackage{color}
\usepackage{bm}
\usepackage{multirow}
\usepackage{enumerate}

\newtheoremstyle{plainNoItalics}{}{}{\normalfont}{}{\bfseries}{.}{ }{}

\theoremstyle{plain}
\newtheorem{thm}{Theorem}[section]

\theoremstyle{plainNoItalics}

\newtheorem{defn}[thm]{Definition}
\newtheorem{rem}[thm]{Remark}
\newtheorem{prop}[thm]{Proposition}
\newtheorem{exa}[thm]{Example}

\newcommand{\f}{\frac}

\newcommand{\beq}{\begin{equation}}
\newcommand{\eeq}{\end{equation}}
\newcommand{\beqa}{\begin{eqnarray}}
\newcommand{\eeqa}{\end{eqnarray}}
\newcommand{\bit}{\begin{itemize}}
\newcommand{\eit}{\end{itemize}}
\newcommand{\bedef}{\begin{defn}}
\newcommand{\edefn}{\end{defn}}
\newcommand{\bpro}{\begin{prop}}
\newcommand{\epro}{\end{prop}}

\newcommand{\Dt}{\Delta t}

\begin{document}

\baselineskip=1.6pc


\begin{center}
{\bf
Parametrized Positivity Preserving Flux Limiters for the High Order Finite Difference WENO Scheme Solving Compressible Euler Equations}
\end{center}
\vspace{.2in}
\centerline{
Tao Xiong \footnote{Department of
Mathematics, University of Houston, Houston, 77204. E-mail:
txiong@math.uh.edu.}
Jing-Mei Qiu \footnote{Department of Mathematics, University of Houston,
Houston, 77204. E-mail: jingqiu@math.uh.edu.
The first and second authors are supported by Air Force Office of Scientific Computing YIP grant FA9550-12-0318, NSF grant DMS-1217008.}
Zhengfu Xu \footnote{Department of
Mathematical Science, Michigan Technological University, Houghton, 49931. E-mail: zhengfux@mtu.edu. Supported by NSF grant DMS-1316662.}
}

\bigskip
\centerline{\bf Abstract}
\vspace{0.2cm}
In this paper, we develop parametrized positivity satisfying flux limiters for the high order finite difference Runge-Kutta weighted essentially non-oscillatory (WENO) scheme solving compressible Euler equations to maintain positive density and pressure. Negative density and pressure, which often leads to simulation blow-ups or nonphysical solutions, emerges from many high resolution computations in some extreme cases.  The methodology we propose in this paper is a nontrivial generalization of the parametrized maximum principle preserving flux limiters for high order finite difference schemes solving scalar hyperbolic conservation laws \cite{mpp_xu, mpp_xuMD, mpp_xqx}. To preserve the maximum principle, the high order flux is limited towards a first order monotone flux, where the limiting procedures are designed by decoupling linear maximum principle constraints. High order schemes with such flux limiters 
are shown to preserve the high order accuracy via local truncation error analysis and by extensive numerical experiments with mild CFL constraints. The parametrized flux limiting approach is generalized to the Euler system to preserve the positivity of density and pressure of numerical solutions via decoupling some nonlinear constraints. Compared with existing high order positivity preserving approaches \cite{zhang2011maximum, zhang2012positivity, zhang2011positivity}, our proposed algorithm is positivity preserving by the design; it is computationally efficient and maintains high order spatial and temporal accuracy in our extensive numerical tests. Numerical tests are performed to demonstrate the efficiency and effectiveness of the proposed new algorithm.

%

\vspace{1cm}
\noindent {{\bf Keywords}: Compressible Euler equations; Positivity preserving; Parametrized flux limiters; High order finite difference method; WENO reconstruction

\newpage

\newpage

\section{Introduction}
\label{sec1}
\setcounter{equation}{0}
\setcounter{figure}{0}
\setcounter{table}{0}

The success of the high order essentially non-oscillatory (ENO) or weight ENO (WENO) methods solving hyperbolic conservation laws has been well documented in the literature \cite{harten1987uniformly, shu1988efficient, liu1996nonoscillatory, Jiang_Shu} and the references therein. At the heart of the high order ENO/WENO schemes solving hyperbolic problem is the robustness, namely stability in the sense of suppressing spurious oscillations around discontinuities. The application of the high order finite difference, finite volume ENO/WENO methods to hyperbolic systems \cite{shu1988efficient, Jiang_Shu}, such as the compressible Euler equations
\begin{eqnarray}
\label{eq:euler}
\left(
\begin{array}{l}
\rho\\
\rho u\\
E
\end{array}
\right)_t+
\left(
\begin{array}{l}
\rho u \\
\rho u + P\\
(E+ P)u
\end{array}
\right)_x=0,
\end{eqnarray}
achieves the goal of suppressing oscillations when discontinuous solution emerges during the time evolution. However, in the extreme case, such as high Mach flow simulation, a slightly different (although equally important) problem is that the high order schemes that we are using might produce solutions with negative density and pressure, which leads to an ill-posed problem, often seen as blow-up of the numerical simulation. The failure of preserving positive density and pressure by the above mentioned schemes in such circumstance pose tremendous difficulty of applying high order schemes to some of the challenging simulations in practice. 

In the earlier work, see \cite{einfeldt1991godunov, linde1997robust, perthame1996positivity} and references included, much attention has been paid to the positivity preservation of schemes up to second order. It wasn't until the recent work by Zhang \& Shu \cite{zhang2010positivity} that arbitrarily high order finite volume WENO and discontinuous Galerkin methods are designed to preserve positivity. The method proposed in \cite{zhang2010positivity} is a successful generalization of their earlier work on the maximum principle preserving (MPP) computations of scalar conservation laws,
see \cite{zhang2011maximum}. Their approach relies on limiting the reconstructed polynomials (finite volume WENO) or representing polynomials (discontinuous Galerkin) around cell averages to be MPP. 
The positivity preserving (PP) finite volume WENO scheme and DG scheme by Zhang \& Shu can be proved to have the designed arbitrary high order accuracy when equipped with proper CFL number. In the later work by the authors \cite{zhang2012positivity}, a PP finite difference WENO method is presented when the density and pressure is strictly greater than a fixed positive constant. In \cite{hu2013positivity}, a flux cut-off limiter method is applied to the high order finite difference WENO method to ensure positive density and pressure. 

In this paper, we continue along the line of research on the parametrized flux limiters proposed in \cite{mpp_xu, mpp_xuMD, mpp_xqx} for high order ENO/WENO methods solving a scalar hyperbolic conservation law
\begin{eqnarray}
\label{hcl} u_t+f(u)_x=0
\end{eqnarray}
subject to the initial condition $u({x}, 0)= u_0 ({x})$.  For this particular family of equations, the solution satisfies a strict maximum principle
\begin{eqnarray}
\label{CMPP}
u_m\le u(x, t) \le u_M \quad \text{if}  \quad u_m\le u_0(x)\le u_M.
\end{eqnarray}
The idea of the parametrized flux limiters for general conservative scheme solving scalar conservation laws is to modify high order numerical fluxes to enforce the discrete maximum principle for the updated solution. In general, a conservative high order scheme with explicit multi-stage Runge-Kutta (RK) time integration for (\ref{hcl}) can be written as
\begin{equation}
\label{Conserative}
u^{n+1}_j=u^{n}_j-\frac{\Delta t}{\Delta x} (\hat H^{rk}_{j+\frac12}-\hat H^{rk}_{j-\frac12}),
\end{equation} 
where $\hat{H}^{rk}_{j\pm\frac12}$ are the corresponding fluxes at the final stage of RK methods.
The MPP properties of high order schemes are realized by taking a convex combination of a high order flux $\hat{H}^{rk}_{j+\frac12}$ and a first order monotone flux $\hat{h}_{j+\frac12}$: $\tilde H^{rk}_{j+\frac12}=\hat{h}_{j+\frac12} + \theta_{j+\frac12} (\hat{H}^{rk}_{j+\frac12}-\hat{h}_{j+\frac12})$, with $\theta_{j+\frac12} \in [0, 1]$. The limiting parameters $\theta_{j+\frac12}$, which measure the change of numerical fluxes, can be found out through decoupling the following MPP constraints that are linear with respect to $\theta_{j\pm\frac12}$,
\begin{equation}
\label{DMP}
u_m\le u^{n+1}_j=u^{n}_j-\frac{\Delta t}{\Delta x} (\tilde H^{rk}_{j+\frac12}-\tilde H^{rk}_{j-\frac12})\le u_M.
\end{equation} 
The similar idea is utilized in this paper in the sense of making sufficient modification of the high order numerical fluxes to ensure that the updated density and pressure are positive.
When such parametrized flux limiters are generalized to preserve the positivity of density and pressure of numerical solutions for Euler equations with source terms, there are several new challenges. 
One of the main difficulties is that the linear MPP constraint (\ref{DMP}) becomes nonlinear for positivity preservation of pressure, which has nonlinear dependence on the density, momentum and energy. We address such challenges by decoupling the nonlinear PP constraint for a `convex set' of the limiting parameters. The proposed approach provides a sufficient condition for preserving positive pressure.  The presence of the source term can also be conveniently handled in the parametrized flux limiting framework. Notice that we only require positivity preservation for the solutions at the final stage of RK method for the sake of preserving the designed high order temporal accuracy. If there are negative density and pressure in intermediate stages of the RK method, the speed of sound is computed by $c=\sqrt{\gamma\frac{|p|}{|\rho|}}$.

Our approach is similar to those very early discussions of the flux limiting approach \cite{boris1973flux, chakravarthy1983high, engquist1980stable, van1974towards, sweby1984high} for the purpose of achieving a total variation diminishing (TVD) property, which is a much stronger stability requirement than the maximum principle. The schemes are expected to be TVD, therefore, most of the schemes are at most of second order accurate.  
To distinguish our work from others' in the context of designing arbitrarily high order schemes, we would like to point out that the method we are proposing only involves the modification of high order numerical fluxes. Another critical difference is that the parametrized flux limiters are only applied to the final stage of the multi-stage RK methods. These new features are designed to produce numerical solutions with positive density and pressure, while allowing for relatively large CFL numbers without sacrificing accuracy in our extensive numerical tests. 
The proposed method is essentially different from those by Zhang \& Shu \cite{zhang2012positivity}, in which the PP property is realized only with fine enough numerical meshes, when the density and the pressure is extremely close to $0$. The flux limiting method we are proposing is also different from the flux cut-off method by Hu \cite{hu2013positivity}, whose approach demands significantly reduced CFL for accuracy as illustrated in their analysis and numerical tests. 
However, the proof of maintaining high order accuracy when the PP flux limiters are applied to the finite difference WENO method solving the Euler system is very difficult. In this paper, we rely on numerical observations to demonstrate the maintenance of high order accuracy. A rigorous proof of that the MPP flux limiters modify the original high order flux with up to third order accuracy for general nonlinear scalar cases is provided in \cite{mpp_xqx} and that with up to fourth order accuracy for linear advection equations is provided in \cite{mpp_vp}.



The paper is organized as follows. In Section \ref{sec2}, we give a brief review of the parametrized MPP flux limiters for high order conservative schemes solving (\ref{hcl}). We then generalize the MPP flux limiters to a scalar problem with source terms. In Section \ref{sec3}, we present the main algorithm of the parametrized PP finite difference WENO RK method for the compressible Euler equation in one and two dimensions. An implementation procedure is given in the presence of source terms. In Section \ref{sec5}, we perform extensive numerical tests to illustrate the effectiveness of the proposed method. We finally conclude in Section \ref{sec6}.

\section{Parametrized MPP flux limiters for scalar equations}
\label{sec2}
\setcounter{equation}{0}
\setcounter{figure}{0}
\setcounter{table}{0}


\subsection{Review of MPP flux limiters for scalar equations}
\label{sec2.1}

For simplicity, we consider a simple one-dimensional hyperbolic conservation equation
\begin{eqnarray}
\label{eq: adv}
u_t+f(u)_x=0, \quad x \in [0, 1],
\end{eqnarray}
with an initial condition $u(x,0) = u_0(x)$ and a periodic boundary condition.
We adopt the following spatial discretization for the domain $[0, 1]$
\[
0 = x_\frac12 < x_\frac32 < \cdots < x_{N+\frac12} = 1,
\]
where $I_j = [x_{j-\frac12}, x_{j+\frac12}]$ has the mesh size $\Delta x = \frac1N$.
Let $u_j(t)$ denote the solution at grid point $x_j = \frac12(x_{j-\frac12}+x_{j+\frac12})$ at continuous time $t$.
The finite difference scheme evolves the point values of the solution in a conservative form
\begin{eqnarray}
\label{eq: semi-discrete}
\frac{d}{dt}u_j(t)+ \frac{1}{\Delta x} (\hat H_{j+1/2}-\hat H_{j-1/2}) = 0.
\end{eqnarray}
The numerical flux $\hat{H}_{j+\frac12}$ in equation \eqref{eq: semi-discrete} can be reconstructed from neighboring flux functions $f(u(x_i, t)),$ $i=j-p, \cdots, j+q$ with high order by WENO reconstructions \cite{Jiang_Shu,shu1998essentially}. By adaptively assigning nonlinear weights to neighboring candidate stencils, the WENO reconstruction preserves high order accuracy of the linear scheme around smooth regions of the solution, while producing a sharp and essentially non-oscillatory capture of discontinuities. Equation \eqref{eq: semi-discrete} can be further discretized in time by a high order time integrator via the method-of-line approach. For example, the scheme with a third order total variation diminishing (TVD) RK time discretization is
\begin{eqnarray}
u_j^{(1)}   &=& u_j^n+\Delta t L(u_j^n),    \nonumber  \\
u_j^{(2)}   &=& u_j^n+\frac{1}{4}\Delta t (L(u_j^n) + L(u_j^{(1)})),  \nonumber \\
u_j^{n+1} &=& u_j^n+\frac{1}{6} \Delta t \bigl( L(u_j^{n})+ L(u_j^{(1)})+4 L(u_j^{(2)})\bigr).
\label{eq:rk3}
\end{eqnarray}
where $u^{(k)}_{j}$ and $u^n_j$ denotes the numerical solution at $x_j$ at $k^{th}$ RK stage and at time $t^n$ respectively.
Let $\Delta t$ be the time step size. $L(u^{(k)}) \doteq -\frac{1}{\Delta x} (\hat H^{(k)}_{j+\f12}-\hat H^{(k)}_{j-\f12})$ with $\hat H^{(k)}_{j+\f12}$ being the numerical flux from 
finite difference WENO reconstruction based on $\{u_j^{(k)}\}_{j=1}^N$ at intermedia RK stages. Equation (\ref{eq:rk3}) in the final stage of RK method can be re-written as
\begin{eqnarray}
\label{eq:rkfinal}
u^{n+1}_j=u^{n}_j-\lambda (\hat H^{rk}_{j+\f12}-\hat H^{rk}_{j-\f12}),
\end{eqnarray}
with $\lambda=\frac{\Delta t}{\Delta x}$ and
\begin{eqnarray}
\hat H^{rk}_{j+\f12} \doteq \frac{1}{6}\left(\hat H^n_{j+\f12}+\hat H^{(1)}_{j+\f12}+4\hat H^{(2)}_{j+\f12}\right).
\label{eq:rkflux}
\end{eqnarray}

The parametrized MPP flux limiters in \cite{mpp_xqx} is based on the finite difference RK WENO scheme for equation \eqref{eq: adv} reviewed above.
Let $u_{m}= \underset{x}{\text{min}}(u(x, 0))$ and $u_{M}=\underset{x}{\text{max}}(u(x, 0))$.
The idea of the parametrized MPP flux limiter is to modify the high order flux $\hat{H}^{rk}_{j\pm\frac12}$ in equation \eqref{eq:rkflux} towards a first order monotone flux denoted as $\hat{h}_{j\pm\frac12}$ by taking
a linear combination of them,
\beq
\label{eq: linear_comb}
\tilde H^{rk}_{j\pm\f12} \doteq \hat{h}_{j\pm\frac12} + \theta_{j\pm\frac12} (\hat{H}^{rk}_{j\pm\frac12}-\hat{h}_{j\pm\frac12}), \quad \theta_{j\pm\frac12} \in [0, 1].
\eeq
the original high order flux $\hat H^{rk}_{j\pm\f12}$ in equation \eqref{eq:rkflux} is then replaced by the modified flux $\tilde H^{rk}_{j\pm\f12}$ above.

To preserve the MPP property, we wish to have $u_{m}\le u^{n+1}_{j} \le u_{M}$ at the final RK stage on each time step, i.e.
\begin{eqnarray}
\label{eq:mpp}
u_{m} \le u^{n}_j-\lambda (\tilde H^{rk}_{j+\f12}-\tilde H^{rk}_{j-\f12}) \le u_{M}.
\end{eqnarray} 
For the parametrized MPP flux limiter, a pair $(\Lambda_{-\f12, {I_j}}, \Lambda_{+\f12, {I_j}})$ needs to be found such that 
any pair $(\theta_{j-\f12}, \theta_{j+\f12}) \in [0, {\Lambda_{-\f12, {I_j}}}]\times [0, {\Lambda_{+\f12, {I_j}}]}$ satisfies (\ref{eq:mpp}). Under such a constraint, $\theta_{j\pm\f12}$ are chosen to be as close to $1$ as possible for accuracy, which is done by the following three steps. Below $\epsilon$ is 
a small positive number to avoid the denominator to be $0$, e.g., $\epsilon=10^{-13}$.
\begin{enumerate}
\item
The right inequality of (\ref{eq:mpp}), that is the maximum value part, can be rewritten as
\begin{eqnarray}
\label {umax}
\lambda \theta_{j-\f12} (\hat H^{rk}_{j-\f12}-\hat h_{j-\f12}) - \lambda \theta_{j+\f12} (\hat H^{rk}_{j+\f12}-\hat h_{j+\f12})-\Gamma^M_j \le 0,
\end{eqnarray}
where $\Gamma^M_j=u_{M}-u_j+\lambda (\hat h_{j+\f12}-\hat h_{j-\f12}) \ge 0$. 
Let $F_{j-\f12}=\hat H^{rk}_{j-\f12}-\hat h_{j-\f12}$, the decoupling of (\ref{umax}) on cell $I_j$ gives:
\begin{enumerate}
\item  If $F_{j-\f12}\le 0$ and $F_{j+\f12}\ge 0$,  let $(\Lambda^M_{-\f12, I_j}, \Lambda^M_{+\f12, I_j})=(1, 1)$.

\item  If $F_{j-\f12}\le 0$ and $F_{j+\f12}  <  0$, let $(\Lambda^M_{-\f12, {I_j}}, \Lambda^M_{+\f12, {I_j}})=(1, \min(1, \frac{\Gamma^M_j}{-\lambda F_{j+\f12}+\epsilon}))$.

\item  If $F_{j-\f12}  > 0$ and $F_{j+\f12}\ge 0$, let $(\Lambda^M_{-\f12, {I_j}}, \Lambda^M_{+\f12, {I_j}})=(\min(1, \frac{\Gamma^M_j}{\lambda F_{j-\f12}+\epsilon}), 1)$.

\item  If $F_{j-\f12}  > 0$ and $F_{j+\f12}  <  0$,
\bit
\item if $(\theta_{j-\f12}, \theta_{j+\f12})=(1, 1)$ satisfies (\ref{umax}), let $(\Lambda^M_{-\f12, {I_j}}, \Lambda^M_{+\f12, {I_j}})=(1, 1)$;
\item otherwise, let
 $(\Lambda^M_{-\f12, {I_j}}, \Lambda^M_{+\f12, {I_j}})=(\frac{\Gamma^M_j}{\lambda F_{j-\f12}- \lambda F_{j+\f12}+\epsilon},\frac{\Gamma^M_j}{\lambda F_{j-\f12}- \lambda F_{j+\f12}+\epsilon} )$.
\eit
\end{enumerate}
\item
The left inequality of (\ref{eq:mpp}), that is the minimum value part, can be rewritten as
\begin{eqnarray}
\label {umin}
0\le \lambda  \theta_{j-\f12} (\hat H^{rk}_{j-\f12}-\hat h_{j-\f12}) - \lambda \theta_{j+\f12} (\hat H^{rk}_{j+\f12}-\hat h_{j+\f12})-\Gamma^m_j,
\end{eqnarray}
where $\Gamma^m_j=u_{m}-u_j+\lambda (\hat h_{j+\f12}-\hat h_{j-\f12}) \le 0$. Similar to the maximum value case, the decoupling of (\ref{umin}) on cell $I_j$ gives:
\begin{enumerate}
\item If $F_{j-\f12}\ge 0$ and $F_{j+\f12}\le 0$, let $(\Lambda^m_{-\f12, I_j}, \Lambda^m_{+\f12, I_j})=(1, 1)$;

\item  If $F_{j-\f12}\ge 0$ and $F_{j+\f12}> 0$, let $(\Lambda^m_{-\f12, {I_j}}, \Lambda^m_{+\f12, {I_j}})=(1, \min(1, \frac{\Gamma^m_j}{-\lambda F_{j+\f12}-\epsilon}))$;

\item  If $F_{j-\f12}< 0$ and $F_{j+\f12}\le 0$, let $(\Lambda^m_{-\f12, {I_j}}, \Lambda^m_{+\f12, {I_j}})=(\min(1, \frac{\Gamma^m_j}{\lambda F_{j-\f12}-\epsilon}), 1)$;

\item  If $F_{j-\f12}< 0$ and $F_{j+\f12}> 0$, 
\bit
\item when $(\theta_{j-\f12}, \theta_{j+\f12})=(1, 1)$ satisfies (\ref{umin}),  
let $(\Lambda^m_{-\f12, {I_j}}, \Lambda^m_{+\f12, {I_j}})=(1, 1)$;
\item otherwise, let $(\Lambda^m_{-\f12, {I_j}}, \Lambda^m_{+\f12, {I_j}})=(\frac{\Gamma^m_j}{\lambda F_{j-\f12}- \lambda F_{j+\f12}-\epsilon},\frac{\Gamma^m_j}{\lambda F_{j-\f12}- \lambda F_{j+\f12}-\epsilon} )$.
\eit
\end{enumerate}
\item 
The locally defined limiting parameter is given as
\begin{eqnarray}
\label{limit1}
\Lambda_{j+\f12}=\min(\Lambda^M_{+\f12, {I_j}}, \Lambda^M_{-\f12, {I_{j+1}}}, \Lambda^m_{+\f12, {I_j}}, \Lambda^m_{-\f12, {I_{j+1}}}), \quad j = 0, \cdots N.
\end{eqnarray}
\end{enumerate}
The flux limiting procedure above guarantees the MPP property of the numerical solution by the design. It is theoretically proved to preserve up to fourth order spatial and temporal accuracy for smooth solutions \cite{mpp_xqx, mpp_vp}. 

\subsection{Scalar advection equations with source terms}
\label{sec2.2}

We consider scalar advection problems with a source term
\begin{eqnarray}
u_t+f(u)_x=s(u).
\label{eq:source1}
\end{eqnarray}
In particular, we consider the class of problems whose solutions enjoy the PP property, that is, the lower bound of the solution is $0$ (such kind of problem might not preserve the MPP property). For example, when $s(u)=-k u$ with a positive $k$, with positive initial values and periodic boundary conditions, the solution satisfies the PP property. 
The flux limiter is designed base on the PP property of a first order scheme
\begin{eqnarray}
\label{eq:1st}
u^{n+1}_j=u^{n}_j-\lambda (\hat h_{j+\f12}-\hat h_{j-\f12})+ \Delta t s(u^n_j),
\end{eqnarray}
under the time step constraint
\begin{equation}
\Delta t \le \frac{\text{CFL } \Delta x}{\lambda_{max}+ s_{max} \Delta x},
\label{eq:timestep}
\end{equation}
where $\lambda_{max}=\max|f'(u)|$ and $s_{max}=\max|s'(u)|$. 

We propose to first modify the source term such that $\tilde{u}^{n+1}_j \ge \epsilon_s$, with $\epsilon_s=\min_j(u^{n+1}_j, 10^{-13})$, where $\{u^{n+1}_j\}$ are positive solutions computed from (\ref{eq:1st}) and $10^{-13}$ is a small positive number related to machine precision. Here $\tilde{u}^{n+1}_j$ is 
\begin{eqnarray}
 \tilde {u}^{n+1}_j=u^{n}_j-\lambda (\hat h_{j+\f12}-\hat h_{j-\f12})+ \Delta t \tilde{s}^{rk}_j,
 \label{eq:source2}
\end{eqnarray}
with
\beq
 \label{eq:mdsource}
\tilde{s}^{rk}_j=r_j (\hat s^{rk}_j-s(u^n_j))+s(u^n_j),
\eeq
and
\begin{eqnarray}
\hat s^{rk}_j \doteq \frac{1}{6}\left(s(u^n_j)+ s(u^{(1)}_j) + 4 s(u^{(2)}_j)\right), 
\label{eq:rksource}
\end{eqnarray}
as in (\ref{eq:rkflux}). $r_j$ is designed by the linear constraints to preserve the PP property of $\{\tilde{u}^{n+1}_j\}_j$. Specifically, 
\[
r_j=
\begin{cases}
\min(\frac{\epsilon_s-u^{n+1}_j}{\Delta t \Delta s_j}, 1),& \quad \text{if } \tilde{\tilde{u}}_j < \epsilon_s \\
1, &\quad \text{otherwise }
\end{cases},
\]
where $\Delta s_j=\hat s^{rk}_j-s(u^n_j)$ and $\tilde{\tilde{u}}_j=u^{n}_j-\lambda (\hat h_{j+\f12}-\hat h_{j-\f12})+ \Delta t \hat s^{rk}_j$.
Next the parametrized MPP flux limiters are applied as in (\ref{eq:mpp}) to satisfy
\begin{eqnarray}
\label{eq:pps}
\epsilon_s \le u^{n}_j-\lambda (\tilde H^{rk}_{j+\f12}-\tilde H^{rk}_{j-\f12})+\Delta t \tilde s^{rk}_j .
\end{eqnarray}
(\ref{eq:pps}) leads to the same decomposed inequality (\ref{umin}) for the minimum value part, only
with $\Gamma^m_j$ given by
\begin{eqnarray}
 \Gamma^m_j&=&\epsilon_s-u_j+\lambda (\hat h_{j+\f12}-\hat h_{j-\f12})-\Delta t \tilde s^{rk}_j \le 0.
\end{eqnarray}
The procedure proposed above for treating equations with a source term is PP by the design, and is shown to maintain high order accuracy by numerical tests in Section~\ref{sec5}.

\section{Parametrized PP flux limiters for compressible Euler equations }
\label{sec3}
\setcounter{equation}{0}
\setcounter{figure}{0}
\setcounter{table}{0}

In this section, we first extend the parametrized MPP flux limiters to PP flux limiters for the compressible Euler equations. We then describe how to generalize
the proposed approach to systems with source terms and to high dimensional systems. In this section, we use letters in bold for vectors.

\subsection{Parametrized positivity preserving flux limiters}
\label{sec3.1}
For compressible Euler equations in one dimension
\begin{eqnarray}
\label{eq:eulers}
{\bf u}_t+{\bf f}({\bf u})_x=0,
\end{eqnarray}
with ${\bf u}=(\rho, \rho u, E)^T$, ${\bf f}({\bf u})=(\rho u, \rho u^2+ p, (E+p)u )^T$,
where $\rho$ is the density, $u$ is the velocity, $p$ is the pressure, $m=\rho u$ is the momentum, $E=\frac{1}{2}\rho u^2+\frac{p}{\gamma-1}$ is the total energy from equation of state (EOS) and $\gamma$
is the ratio of specific heat ($\gamma=1.4$ for the air).
Denote $\hat {\bf h}_{j+\f12}$ to be a first order monotone flux, and $\hat {\bf H}^{rk}_{j+\f12}$ to be the linear combinations of fluxes from multiple RK stages, similar to equation (\ref{eq:rkflux}), but in a component-by-component fashion.
For positivity preserving, we are seeking the flux limiters of the type
\begin{eqnarray}
\label{mhrk}
\tilde{\bf H}^{rk}_{j+\f12}=\theta_{j+\f12} (\hat{\bf H}^{rk}_{j+\f12}-\hat {\bf h}_{j+\f12})+\hat {\bf h}_{j+\f12}
\end{eqnarray}
such that 
\begin{eqnarray}
\label{eq:pp}
\begin{cases}
\rho^{n+1}_{j}>0,\\
 p^{n+1}_j>0,
\end{cases}
\end{eqnarray}
for the updated solution 
\begin{eqnarray}
\label{eq:rkeuler}
{\bf u}^{n+1}_j={\bf u}^{n}_j-\lambda (\tilde{\bf H}^{rk}_{j+\f12}-\tilde{\bf H}^{rk}_{j-\f12}).
\end{eqnarray}
In the parametrized flux limiters' framework, 
a pair of $(\Lambda_{-\f12, {I_j}}, \Lambda_{+\f12, {I_j}})$ is found such that 
the updated solution satisfies (\ref{eq:pp}) for any $(\theta_{j-\f12}, \theta_{j+\f12}) \in [0, {\Lambda_{-\f12, {I_j}}}]\times [0, {\Lambda_{+\f12, {I_j}}]}$. 
The high order flux $\hat {\bf H}^{rk}_{j+\f12}$ is modified by (\ref{mhrk}) to preserve positive density and pressure. In simulations, preserving positivity is implemented by 
\begin{eqnarray}
\label{eq:pp2}
\begin{cases}
\rho^{n+1}_{j}\ge\epsilon_{\rho},\\
 p^{n+1}_j\ge\epsilon_p.
\end{cases}
\end{eqnarray}
where we introduce small positive numbers $\epsilon_{\rho}$ defined by $\min_{j}(\rho^{n+1}_j, 10^{-13})$ and $\epsilon_p$} defined by $\min_{j}(p^{n+1}_j, 10^{-13})$. $\rho^{n+1}_j$ and $p^{n+1}_j$ are positive density and pressure obtained by the first order monotone scheme and $10^{-13}$ is related to the machine precision.  
Let us denote the first order monotone flux by $\hat {\bf h}({\bf u})=(f^\rho, f^m, f^E)^T$, similarly $\hat {\bf H}^{rk}=(\hat f^\rho, \hat f^m, \hat f^E)^T$ and
$\tilde{\bf H}^{rk}=(\tilde f^{\rho}, \tilde f^{m}, \tilde f^{E})^T$. 
The proposed process can be dissected into two steps. 
%
\begin{enumerate}
\item Find the limiting parameters $\theta_{j\pm\f12}$ to preserve the positivity of the density,
\begin{eqnarray}
\label{density}
\rho^{n+1}_j=\rho^n_j-\lambda (\tilde f^{\rho}_{j+\f12}-\tilde f^{\rho}_{j-\f12}).
\end{eqnarray}
Thus, the limiting parameters $\theta_{j\pm\f12}$ are found to satisfy 
\begin{eqnarray}
\label{d1}
\epsilon_{\rho} \le \Gamma_j-\lambda (\theta_{j+\f12} (\hat f^{\rho}_{j+\f12}-f^{\rho}_{j+\f12}) -\theta_{ j-\f12} (\hat f^{\rho}_{j-\f12}-f^{\rho}_{j-\f12})),
\end{eqnarray}
which is equivalent to 
\begin{eqnarray}
\label{d2}
 0 \le \Gamma_j-\epsilon_{\rho}-\lambda (\theta_{j+\f12} (\hat f^{\rho}_{j+\f12}-f^{\rho}_{j+\f12}) -\theta_{ j-\f12} (\hat f^{\rho}_{j-\f12}-f^{\rho}_{j-\f12})),
\end{eqnarray}
where $\Gamma_j=\rho^n_j-\lambda (f^{\rho}_{j+\f12}- f^{\rho}_{j-\f12}) \ge \epsilon_{\rho}$. 
A pair of limiting parameters $(\Lambda^{\rho}_{-\f12, {I_j}}, \Lambda^{\rho}_{+\f12, {I_j}})$ for the positive density of (\ref{d2}) can be identified by a similar procedure as described in Section~\ref{sec2.1}. We can define a set for the positive density $\rho^{n+1}_j$
\begin{eqnarray}
\label{Srho}
S_{\rho} =\{ (\theta_{j-\f12}, \theta_{j+\f12}): 0\le \theta_{j-\f12} \le \Lambda^{\rho}_{-\f12, {I_j}}, 0\le \theta_{j+\f12} \le   \Lambda^{\rho}_{+\f12, {I_j}} \},
\end{eqnarray}
which is plotted as the rectangle bounded by the dash line in Figure~\ref{set}. 
\item Find the limiting parameters $\theta_{j\pm\f12}$ within the region $S_{\rho}$ to preserve the positivity of the pressure.
We seek a sufficient condition such that the pressure given by (\ref{eq:rkeuler}) satisfies
\begin{eqnarray}
\label{pre}
p^{n+1}_j (\theta_{j-\f12}, \theta_{j+\f12})=(\gamma-1)\left(E^{n+1}_j-\frac{1}{2} \frac{(m^{n+1}_j )^2}{\rho^{n+1}_j}\right)\ge \epsilon_p.
\end{eqnarray}
The decoupling of (\ref{pre}) for ($\theta_{j-\f12}, \theta_{j+\f12}$) is different from 
the scalar case since the principal variables are nonlinearly dependent on each other. However the idea is 
still to separate $\theta_{j-\f12}$ and $\theta_{j+\f12}$. Since $\rho^{n+1}_j \ge \epsilon_{\rho}$ is guaranteed by the previous step, we first put the concave property of pressure \cite{zhang2010positivity} in the following remark for future reference:
\begin{rem}
\label{conv}
The pressure as a function of $(\rho , m, E)$ is concave, i.e., $p(\alpha {\bf U_1}+(1-\alpha) {\bf U_2})\ge \alpha p({\bf U_1})+(1-\alpha) p({\bf U_2})$ for $0\le \alpha \le 1$
if $\rho_1, \rho_2> 0$. Therefore $p^{n+1}_j (\theta_{j-\f12}, \theta_{j+\f12})$ is a concave function of $(\theta_{j-\f12}, \theta_{j+\f12})$ on $S_\rho$ due to the linear dependence of $(\rho^{n+1}_j , m^{n+1}_j, E^{n+1}_j)$ on $(\theta_{j-\f12}, \theta_{j+\f12})$. Therefore, if 
$p^{n+1}_j (\vec{\theta}^l) \ge \epsilon_p$, with $\vec{\theta}^l = (\theta^l_{j-\f12}, \theta^l_{j+\f12})$ for $l=1, 2$,
then 
$
p^{n+1}_j (\alpha \vec{\theta}^1 + (1-\alpha) \vec{\theta}^2) \ge \epsilon_p, \quad 0\le\alpha\le1.
$

\end{rem}
We define an admissible set
\begin{eqnarray}
\label{pAD}
S_\theta =\{(\theta_{j-\f12}, \theta_{j+\f12})\in S_\rho: (\theta_{j-\f12}, \theta_{j+\f12}) \text{ satisfies } (\ref{pre})\}.
\end{eqnarray}
$S_\theta$ is a convex set thanks to Remark \ref{conv}. Let the three vertices of the rectangle $S_\rho$ other than $(0, 0)$ be denoted by
\begin{eqnarray}
\label{vert}
A^1=(0,  \Lambda^{\rho}_{+\f12, {I_j}}),\quad
A^2=(\Lambda^{\rho}_{-\f12, {I_j}}, 0),\quad
A^3=(\Lambda^{\rho}_{-\f12, {I_j}} ,  \Lambda^{\rho}_{+\f12, {I_j}}),
\end{eqnarray}
see Figure~\ref{set}. Based on the concave property in Remark~\ref{conv}, we propose the following way of decoupling (\ref{pre}). 
\begin{enumerate}
\item For i=1, 2, 3, if $p(A^i)\ge \epsilon_p$, let $B^i =A^i$; otherwise find $r$ such that $p(r A^i)\ge \epsilon_p$ and let $B^i= r A^i$. The three $B^i$'s and $(0, 0)$ form a convex polygonal region, denoted as $S_p$,  inside $S_\theta$. Such convex polygonal region $S_p$ is outlined by the dash dot line in Figure~\ref{set}. 
\item We define the decoupling rectangle, as a subset of $S_p$, to be
\begin{eqnarray}
\label{dreg}
R_{\rho, p}=[0, \min(B^2_1, B^3_1)]\times [0, \min(B^1_2, B^3_2)],
\end{eqnarray}
see the region outlined by the solid line in Figure~\ref{set}.
That is, within $S_p$, we find the decoupling rectangle $R_{\rho, p}$ with left-bottom node on $(0, 0)$ and right-top node $(\Lambda_{-\f12, {I_j}}, \Lambda_{+\f12, {I_j}})$ as close to $(1, 1)$ as possible to best preserve the accuracy while achieving the PP property of high order numerical schemes. 
Let
\begin{eqnarray}
\label{dcp}
(\Lambda_{-\f12, {I_j}}, \Lambda_{+\f12, {I_j}})=(\min(B^2_1, B^3_1), \min(B^1_2, B^3_2)).
\end{eqnarray}
\end{enumerate}
\end{enumerate}
Finally, similar to equation \eqref{limit1} for the MPP flux limiters, the locally defined limiting parameter is given as 
$\theta_{j+\f12}=\min(\Lambda_{-\f12,{I_j}}, \Lambda_{+\f12,{I_{j+1}}})$.

\begin{figure}
\centering
\includegraphics[totalheight=3.5in]{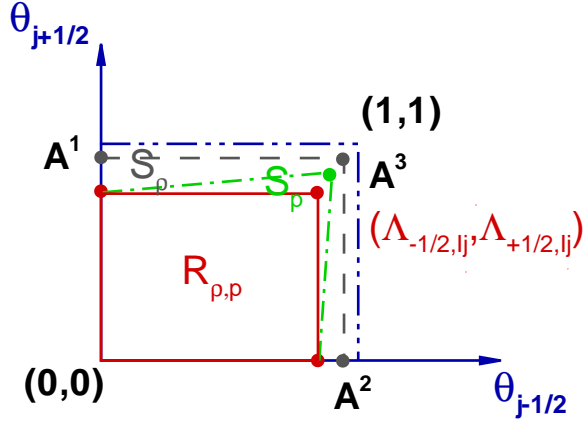}
\caption{The decoupling rectangle $R_{\rho,p}$ (bounded by the solid line) with the right-top node $(\Lambda_{-\f12, {I_j}}, \Lambda_{+\f12, {I_j}})$. $S_{\rho}$ is the rectangle bounded by the dash line. $S_p$ is the polygonal bounded by the dash dot line.}
\label{set}
\end{figure}

\begin{rem}
\label{geos}
The limiter above can preserve positive density and pressure by its design due to the two sufficient conditions (\ref{d1}) and (\ref{pre}).
For general equation of state, if $\rho> 0$, then $p>0 \Leftrightarrow e>0$, where the internal energy $e$ can always be written as a concave function of $(\rho, m, E)^T$ similarly as (\ref{pre}) \cite{zhang2012positivity}. Similar procedure can be followed for PP property of numerical solutions.
\end{rem}

\subsection{Extension to Euler system with source term}
\label{sec3.2}
The compressible Euler equations may come with source terms in the form of
\begin{eqnarray}
\label{eq:eulersource}
{\bf u}_t+{\bf f}({\bf u})_x={\bf s}({\bf u}),
\end{eqnarray}
For example, four kinds of source terms were discussed in \cite{zhang2011positivity}: geometric, gravity, chemical
reaction and radiative cooling. The PP flux limiters can be applied by the following three steps.
\begin{enumerate}
\item Choose a time step, such that the first order scheme (\ref{step1}) is PP,
\begin{eqnarray}
{\bf u}^{n+1}_j={\bf u}^{n}_j-\lambda (\hat {\bf h}_{j+\f12}-\hat {\bf h}_{j-\f12})+\Dt {\bf s}({\bf u}^n_j).
\label{step1}
\end{eqnarray}
\item Find $r$ such that the scheme (\ref{step2}) with the modified source terms is PP
\begin{eqnarray}
{\bf u}^{n+1}_j={\bf u}^{n}_j-\lambda (\hat {\bf h}_{j+\f12}-\hat {\bf h}_{j-\f12})+\Dt \tilde {\bf s}^{rk}_j,
\label{step2}
\end{eqnarray}
with $\tilde {\bf s}^{rk}_j = r(\hat {\bf s}^{rk}_j-{\bf s}({\bf u}^n_j))+{\bf s}({\bf u}^n_j)$, $\hat {\bf s}^{rk}_j$ is similarly defined as (\ref{eq:rksource}) component-by-component.
\item Finally find $\theta_{j\pm\f12}$ for the modified high order flux $\tilde{\bf H}^{rk}_{j+\f12}$, such that (\ref{step3}) is PP
\begin{eqnarray}
{\bf u}^{n+1}_j={\bf u}^{n}_j-\lambda (\tilde {\bf H}^{rk}_{j+\f12}-\tilde {\bf H}^{rk}_{j-\f12})+\Dt \tilde {\bf s}^{rk}_j.
\label{step3}
\end{eqnarray}
The procedure is similar as in the previous subsection.
\end{enumerate}

\subsection{Extension to the multi-dimensional Euler system}
\label{sec3.3}
In this subsection, we extend the previously proposed PP flux limiters to Euler equations 
in two-dimensions
\begin{equation}
 {\bf u}_t+{\bf f}({\bf u})_x+{\bf g}({\bf u})_y=0,
 \label{eq:euler2d}
\end{equation}
with ${\bf u}=(\rho, m_u, m_v, E)^T$, ${\bf f}({\bf u})=(m_u, \rho u^2+p, \rho u v, (E+p)u)^T$ and
${\bf g}({\bf u})=(m_v, \rho u v, \rho v^2+p, (E+p)v)^T$. 
$\rho$ is the density, $u$ is the velocity in $x$ direction, $v$ is the velocity in $y$ direction,
$p$ is the pressure, $m_u=\rho u$ and $m_v=\rho v$ are the momenta, $E=\frac{1}{2}\rho u^2+\frac{1}{2}\rho v^2+\frac{p}{\gamma-1}$ is the total energy and $\gamma$ is the ratio of specific heat. 

The high order finite difference scheme with PP flux limiters at the final stage of a RK time discretization
is given by
\begin{eqnarray}
 {\bf u}^{n+1}_{i,j}={\bf u}^{n}_{i,j}-\lambda_x(\tilde {\bf H}^{rk}_{i+\f12,j}-\tilde {\bf H}^{rk}_{i-\f12,j})
 -\lambda_y(\tilde {\bf G}^{rk}_{i,j+\f12}- \tilde {\bf G}^{rk}_{i,j-\f12}),
\end{eqnarray}
with
\begin{eqnarray}
\label{mhrk2d}
 \tilde {\bf H}^{rk}_{i+\f12,j}&=&\theta_{i+\f12,j}(\hat {\bf H}^{rk}_{i+\f12,j}-\hat {\bf h}_{i+\f12,j})+\hat {\bf h}_{i+\f12,j}, \\
 \label{mgrk2d}
 \tilde {\bf G}^{rk}_{i,j+\f12}&=&\theta_{i,j+\f12}(\hat {\bf G}^{rk}_{i,j+\f12}-\hat {\bf g}_{i,j+\f12})+\hat {\bf g}_{i,j+\f12}, 
\end{eqnarray}
where $\hat {\bf H}^{rk}_{i+\f12,j}$ and $\hat {\bf G}^{rk}_{i,j+\f12}$ are linear combinations of fluxes from multiple RK stages similarly as (\ref{eq:rkflux}) in the scalar case but in a component-wise fashion, $\hat {\bf h}_{i+\f12,j}$ and $\hat {\bf g}_{i,j+\f12}$ are first order monotone fluxes.

Similar to the 1D case, we find the four parametrized limiters $\Lambda^{\rho}_{L, {I_{ij}}}$, $\Lambda^{\rho}_{R,{I_{ij}}}$,
$\Lambda^{\rho}_{U, {I_{ij}}}$ and $\Lambda^{\rho}_{D,{I_{ij}}}$, such that for all $\theta_{i\pm\f12,j}$ and $\theta_{i,j\pm\f12}$
in the set
\begin{align}
\label{Srho2d}
S_{\rho} = &\{ (\theta_{i-\f12,j}, \theta_{i+\f12,j},\theta_{i,j-\f12}, \theta_{i,j+\f12}): 0\le \theta_{i-\f12,j} \le  \Lambda^{\rho}_{L, {I_{ij}}}, \nonumber \\
 &0\le \theta_{i+\f12,j} \le   \Lambda^{\rho}_{R, {I_{ij}}},0\le \theta_{i,j-\f12} \le  \Lambda^{\rho}_{D, {I_{ij}}},
0\le \theta_{i,j+\f12} \le   \Lambda^{\rho}_{U, {I_{ij}}} \}
\end{align}
we have $\rho^{n+1}_{i,j}\ge\epsilon_\rho$. With the positive density $\rho^{n+1}_{i,j}$, the pressure is updated by the constraint
\begin{align}
\label{pre2d}
p^{n+1}_{i,j}& (\theta_{i-\f12,j}, \theta_{i+\f12,j},\theta_{i,j-\f12}, \theta_{i,j+\f12})=\nonumber \\
&(\gamma-1)(E^{n+1}_{i,j} -\frac{1}{2} \frac{((m_u)^{n+1}_{i,j} )^2+((m_v)^{n+1}_{i,j} )^2}{\rho^{n+1}_{i,j}})\ge \epsilon_p.
\end{align}
Let the convex admissible set for positive pressure be
\begin{align}
\label{pAD2d}
S_\theta =\{(\theta_{i-\f12,j}, \theta_{i+\f12,j},\theta_{i,j-\f12}, \theta_{i,j+\f12})\in S_\rho: 
(\theta_{i-\f12,j}, \theta_{i+\f12,j},\theta_{i,j-\f12}, \theta_{i,j+\f12}) \text{ satisfies } (\ref{pre2d})\}
\end{align}
Let the sixteen vertices of $S_\rho$ denoted by
\begin{equation}
 A^{k_1,k_2,k_3, k_4}=(k_1 \Lambda^{\rho}_{L,{I_{ij}}}, k_2 \Lambda^{\rho}_{R,{I_{ij}}},k_3 \Lambda^{\rho}_{D,{I_{ij}}},k_4 \Lambda^{\rho}_{U,{I_{ij}}}),
\end{equation}
with $k_1, k_2, k_3, k_4$ to be $0$ or $1$. We decouple (\ref{pre2d}) in the following way:
\begin{enumerate}
\item For $(k_1,k_2,k_3,k_4)\neq(0,0,0,0)$, if $p(A^{k_1,k_2,k_3,k_4})\ge \epsilon_p$, let $B^{k_1,k_2,k_3,k_4} =A^{k_1,k_2,k_3,k_4}$;
otherwise find $r$ such that $P(r A^{k_1,k_2,k_3,k_4})\ge \epsilon_p$ and let $B^{k_1,k_2,k_3,k_4} =r A^{k_1,k_2,k_3,k_4}$.
The 15 $B^{k_1,k_2,k_3,k_4}$'s with the origin $(0, 0,0,0)$ form a four dimensional polyhedra inside $S_{\theta}$;
\item The decoupling tesseract can be defined by
\begin{align}
\label{dreg2d}
R_{\rho, p}=&[0, \min(B^{1,1,1,0}_1, B^{1,1,0,1}_1,B^{1,0,1,1}_1)] \times [0, \min(B^{1,1,1,0}_2, B^{1,1,0,1}_2,B^{0,1,1,1}_2)] \nonumber \\
\times&  [0, \min(B^{1,1,1,0}_3, B^{1,0,1,1}_3, B^{0,1,1,1}_3)]\times [0, \min(B^{1,1,0,1}_4, B^{1,0,1,1}_4,B^{0,1,1,1}_4)].
\end{align}
Let
\begin{align}
\label{dcp2d}
(\Lambda_{L, {I_{ij}}}&, \Lambda_{R, {I_{ij}}},\Lambda_{D, {I_{ij}}}, \Lambda_{U, {I_{ij}}}) =
(\min(B^{1,1,1,0}_1, B^{1,1,0,1}_1,B^{1,0,1,1}_1), \min(B^{1,1,1,0}_2, B^{1,1,0,1}_2, \nonumber \\  
&B^{0,1,1,1}_2),\min(B^{1,1,1,0}_3, B^{1,0,1,1}_3, B^{0,1,1,1}_3), \min(B^{1,1,0,1}_4, B^{1,0,1,1}_4,B^{0,1,1,1}_4)).
\end{align}
\end{enumerate}
Finally, similar to equation \eqref{limit1} for the MPP flux limiters, the locally defined limiting parameter is given as 
$\theta_{i+\f12,j}=\min(\Lambda_{L,{I_{ij}}}, \Lambda_{R,{I_{i+1,j}}})$ and
$\theta_{i,j+\f12}=\min(\Lambda_{D,{I_{ij}}}, \Lambda_{U,{I_{i,j+1}}})$.
\begin{rem}
For two dimensional compressible Euler equations with source terms, it can be done similarly as the one dimensional case.
\end{rem}
\section{Numerical simulations}
\label{sec5}
\setcounter{equation}{0}
\setcounter{figure}{0}
\setcounter{table}{0}

In this section, we will use the 5th order finite difference WENO scheme for space discretization \cite{Jiang_Shu} and a 4th order Runge-Kutta time discretization \cite{shu1988efficient}, denote as ``WENO5RK4'', 
with the proposed PP flux limiters for simulating the compressible Euler equations. 
Here a 4th order RK time discretization is adopted for better observation of accuracy
by taking the time step to be $\Delta t=\text{CFL } \Delta x$. Most of the tests are from 
\cite{zhang2012positivity}. Below, $\text{CFL }=0.6$ unless otherwise specified.

\begin{exa}(Accuracy test for a scalar problem with a source term.)
\label{ex2}
We consider $u_t+u_x=-u$ with the initial condition 
\begin{equation*}
u(x,0)=\sin^4(x),
\end{equation*}
and the periodic boundary condition. The exact solution is given by
\begin{equation*}
u(x,t)=e^{-t}\sin^4(x-t).
\end{equation*}
The minimum value of the exact solution is $u_{m}=0$. This example is used to test the PP property
and accuracy of dealing with a source term. In Table \ref{tab2}, we can see the PP property 
is preserved and the 5th order accuracy has been maintained.

\begin{table}
\centering
\caption{Example \ref{ex2}. A scalar advection problem with a source term at $T=0.1$. $v_{min}$ is the minimum value of the numerical solution.}
\vspace{0.2cm}
  \begin{tabular}{|c|c|c|c|c|c|c|}
    \hline
     & N  & $L^1$ error &    order   & $L^\infty$ error & order & $v_{min}$ \\\hline
\multirow{5}{*}{without limiters}
  &  40 &     3.36E-04 &       --&     8.78E-04 &       --&    -1.35E-04 \\ \cline{2-7}
  &  80 &     2.03E-05 &     4.05&     1.24E-04 &     2.82&    -1.05E-05 \\ \cline{2-7}
  & 160 &     6.75E-07 &     4.91&     6.25E-06 &     4.31&    -1.88E-06 \\ \cline{2-7}
  & 320 &     1.67E-08 &     5.34&     1.29E-07 &     5.60&    -3.02E-09 \\ \cline{2-7}
  & 640 &     4.30E-10 &     5.28&     2.60E-09 &     5.63&    -5.28E-11 \\ \hline
\multirow{5}{*}{with limiters}
  &  40 &     3.25E-04 &       --&     8.66E-04 &       --&     5.67E-15 \\ \cline{2-7}
  &  80 &     1.92E-05 &     4.08&     1.17E-04 &     2.89&     1.18E-05 \\ \cline{2-7}
  & 160 &     6.38E-07 &     4.91&     6.25E-06 &     4.22&     3.01E-16 \\ \cline{2-7}
  & 320 &     1.67E-08 &     5.26&     1.29E-07 &     5.60&     6.33E-10 \\ \cline{2-7}
  & 640 &     4.31E-10 &     5.28&     2.60E-09 &     5.63&     3.46E-11 \\ \hline
  \end{tabular}
\label{tab2}
\end{table}
\end{exa}

\begin{exa}(Accuracy test for the global Lax-Friedrichs flux.) 
\label{ex22}
We consider the Burgers' equation with the initial condition
\begin{equation*}
u(x,0)=(1+\sin(x))/2
\end{equation*}
and a periodic boundary condition. 
We consider the WENO5RK4 scheme with the global Lax-Friedrichs (LxF) fluxes. Let
\beq
\label{eq: g_LxF}
f^\pm_i=\f12(f(u^n_i)\pm \alpha u^n_i),\quad i=j-p, \cdots, j+q,
\eeq
with $\alpha\ge \max_{u} |f'(u)|$. The numerical flux $\hat H_{j+\f12}=f^-_{j+\f12}+f^+_{j+\f12}$ in (\ref{eq: semi-discrete}), where $f^\pm_{j+\f12}$ are reconstructed based on WENO schemes from (\ref{eq: g_LxF}) with the corresponding upwind mechanism. We numerically investigate the time step restriction for maintaining high order accuracy using the global Lax-Friedrichs flux, since it is frequently used in the computation of the Euler system. In \cite{mpp_xqx}, local truncation analysis is performed to prove that MPP flux limiters can maintain up to third order accuracy of the original scheme with no additional CFL constraint (i.e. $\text{CFL}\le1$) when the upwind flux is used. However, when the global LxF flux with extra large $\alpha$ in equation \eqref{eq: g_LxF} is used, there is a mild time step restriction with $\text{CFL}\le0.886$. It is technically challenging to theoretically estimate such time step restriction for maintaining high order accuracy (e.g. fifth order) of the MPP flux limiters even for scalar equations, therefore we rely on extensive numerical tests.

We consider the scheme with the global Lax-Friedrichs flux with extra large $\alpha = 1.3$ (greater than $\max_u|f'(u)|=1$). The time step is chosen to be $\Delta t=\text{CFL} \Delta x / \alpha$. In Table \ref{tab22}, we show that for the 5th order linear scheme (linear weights instead of nonlinear weights in WENO5) with the 4th order Runge-Kutta time discretization, when $\text{CFL}=0.886$, the 5th order accuracy is maintained with the MPP flux limiters. In fact, $\text{CFL}=0.886$ works for all other $\alpha$'s we tested, 
the results are not listed here to save space.

\begin{table}
\centering
\caption{Example \ref{ex22}. Burgers' equation at $T=0.2$. $\alpha=1.3$ for the global LxF flux (\ref{eq: g_LxF}). $\Delta t=0.886 \Delta x / \alpha $. $u_{max}-v_{max}$ is the difference of the maximum values between the numerical solution and the exact solution.}
\vspace{0.2cm}
  \begin{tabular}{|c|c|c|c|c|c|c|}
    \hline
     & N  & $L^1$ error &    order   & $L^\infty$ error & order & $u_{max}-v_{max}$ \\\hline
\multirow{5}{*}{without limiters}
  &   40 &     2.05E-04 &       --&     2.76E-03 &       --  &  -5.47E-06  \\  \cline{2-7}
  &   80 &     1.20E-05 &     4.09&     2.33E-04 &     3.57  &  -2.24E-07  \\  \cline{2-7}
  &  160 &     4.32E-07 &     4.79&     9.68E-06 &     4.59  &  -1.98E-08  \\  \cline{2-7}
  &  320 &     1.38E-08 &     4.97&     3.15E-07 &     4.94  &  -1.37E-09  \\  \cline{2-7}
  &  640 &     4.39E-10 &     4.98&     1.01E-08 &     4.97  &  -8.92E-11  \\  \hline
\multirow{5}{*}{with limiters}
  &   40 &     2.06E-04 &       --&     2.76E-03 &       --  &   7.33E-06  \\  \cline{2-7}
  &   80 &     1.20E-05 &     4.10&     2.33E-04 &     3.57  &   9.99E-14  \\  \cline{2-7}
  &  160 &     4.32E-07 &     4.79&     9.68E-06 &     4.59  &   1.00E-13  \\  \cline{2-7}
  &  320 &     1.38E-08 &     4.97&     3.15E-07 &     4.94  &   1.00E-13  \\  \cline{2-7}
  &  640 &     4.39E-10 &     4.98&     1.01E-08 &     4.97  &   9.99E-14  \\  \hline
  \end{tabular}
\label{tab22}
\end{table}
\end{exa}

\begin{exa}{(Accuracy test for 2D vortex evolution problem.)}
\label{ex1}
We consider the vortex evolution problem \cite{hu1999weighted} to test the accuracy.
For this problem, the mean flow is $\rho=p=u=v=1$ and is added by an isentropic vortex perturbation centered at $(x_0,y_0)$ in $(u,v)$ with $T=p/\rho$, no perturbation in entropy $S=p/\rho^\gamma$,
\begin{eqnarray}
 (\delta u,\delta v)=\frac{\varepsilon_{vortex}}{2\pi}e^{0.5(1-r^2)}(-\bar{y},\bar{x}),
 \quad \delta T=-\frac{(\gamma-1)\epsilon^2}{8\gamma\pi^2}e^{(1-r^2)},
 \quad \delta S=0,
\end{eqnarray}
where $(\bar{x},\bar{y})=(x-x_0,y-y_0)$, $r^2=\bar{x}^2+\bar{y}^2$.

The computational domain is taken to be $[-5, 15]\times[-5, 15]$ and $(x_0,y_0)=(5,5)$.
The boundary condition is periodic. $\gamma=1.4$ and the vortex strength is $\varepsilon_{vortex}=10.0828$ as in \cite{zhang2012positivity}. The exact solution is the passive convection of the vortex with the mean flow. The lowest density and pressure of the exact solution are $7.8\times 10^{-15}$ and $1.7\times 10^{-20}$.

$\epsilon_{WENO}$ in the nonlinear WENO weights is chosen to be $10^{-5}$, which is between $10^{-2}$ and $10^{-6}$ \cite{hu1999weighted}. In Table \ref{tab1}, we can clearly observe the 5th order accuracy with the PP flux limiters. 


\begin{table}
\centering
\caption{Example \ref{ex1}. Vortex evolution problem at $T=0.01$. $\epsilon_{WENO}=10^{-5}$. $\rho_{min}$ and
$p_{min}$ are the minimum density and pressure of the numerical solution respectively.}
\vspace{0.2cm}
  \begin{tabular}{|c|c|c|c|c|c|c|c|}
    \hline
    &  N  & $L^1$ error &    order   & $L^\infty$ error & order & $\rho_{min}$ & $p_{min}$\\ \hline
\multirow{5}{*}{without limiters}
 &   64 &     1.49E-04 &       -- &     5.25E-02 &       -- &     -9.10E-05 &      2.79E-04  \\ \cline{2-8}
 &  128 &     1.57E-06 &     6.57 &     5.39E-04 &     6.61 &      3.04E-06 &     -2.16E-06  \\ \cline{2-8}
 &  256 &     1.29E-07 &     3.60 &     1.37E-04 &     1.97 &     -2.83E-06 &     -4.68E-07  \\ \cline{2-8}
 &  512 &     4.69E-09 &     4.79 &     3.37E-06 &     5.35 &      2.42E-07 &      1.27E-08  \\ \cline{2-8}
 & 1024 &     1.15E-10 &     5.35 &     7.92E-08 &     5.41 &      1.31E-08 &      1.87E-10  \\ \hline
\multirow{5}{*}{with limiters}
 &   64 &     1.49E-04 &       -- &     5.25E-02 &       -- &      6.30E-04 &      1.91E-04  \\ \cline{2-8}
 &  128 &     1.57E-06 &     6.57 &     5.39E-04 &     6.61 &      3.72E-06 &      1.00E-13  \\ \cline{2-8}
 &  256 &     1.32E-07 &     3.57 &     1.30E-04 &     2.05 &      2.42E-07 &      5.36E-10  \\ \cline{2-8}
 &  512 &     4.69E-09 &     4.81 &     3.37E-06 &     5.27 &      2.42E-07 &      1.27E-08  \\ \cline{2-8}
 & 1024 &     1.15E-10 &     5.35 &     7.92E-08 &     5.41 &      1.31E-08 &      1.87E-10  \\ \hline

  \end{tabular}
\label{tab1}
\end{table}

\end{exa}

\begin{exa}{1D low density and low pressure problems.}
\label{ex3}
We consider two 1D low density and low pressure problems for the ideal gas.
The first one is a 1D Riemann problem, the initial condition is $\rho_L=\rho_R=7$,
$u_L=-1$, $u_R=1$, $p_L=p_R=0.2$ and $\gamma=1.4$, which is a double rarefaction problem. 
The exact solution contains vacuum. In Fig. \ref{fig52} (left), we show the results with 
the PP flux limiters at $T=0.6$ on a mesh size of $\Delta x=1/200$.

The second one is the 1D Sedov blast wave. For the initial condition, the density is $1$, 
the velocity is $0$, the total energy is $10^{-12}$ everywhere except in the center cell, which 
is a constant $E_0/\Delta x$ with $E_0=3200000$. $\gamma=1.4$. In Fig. \ref{fig52} (right),
we show the results with the PP flux limiters at $T=0.001$ on a mesh size of $\Delta x=1/200$.

\begin{figure}
\centering
\includegraphics[totalheight=2.0in]{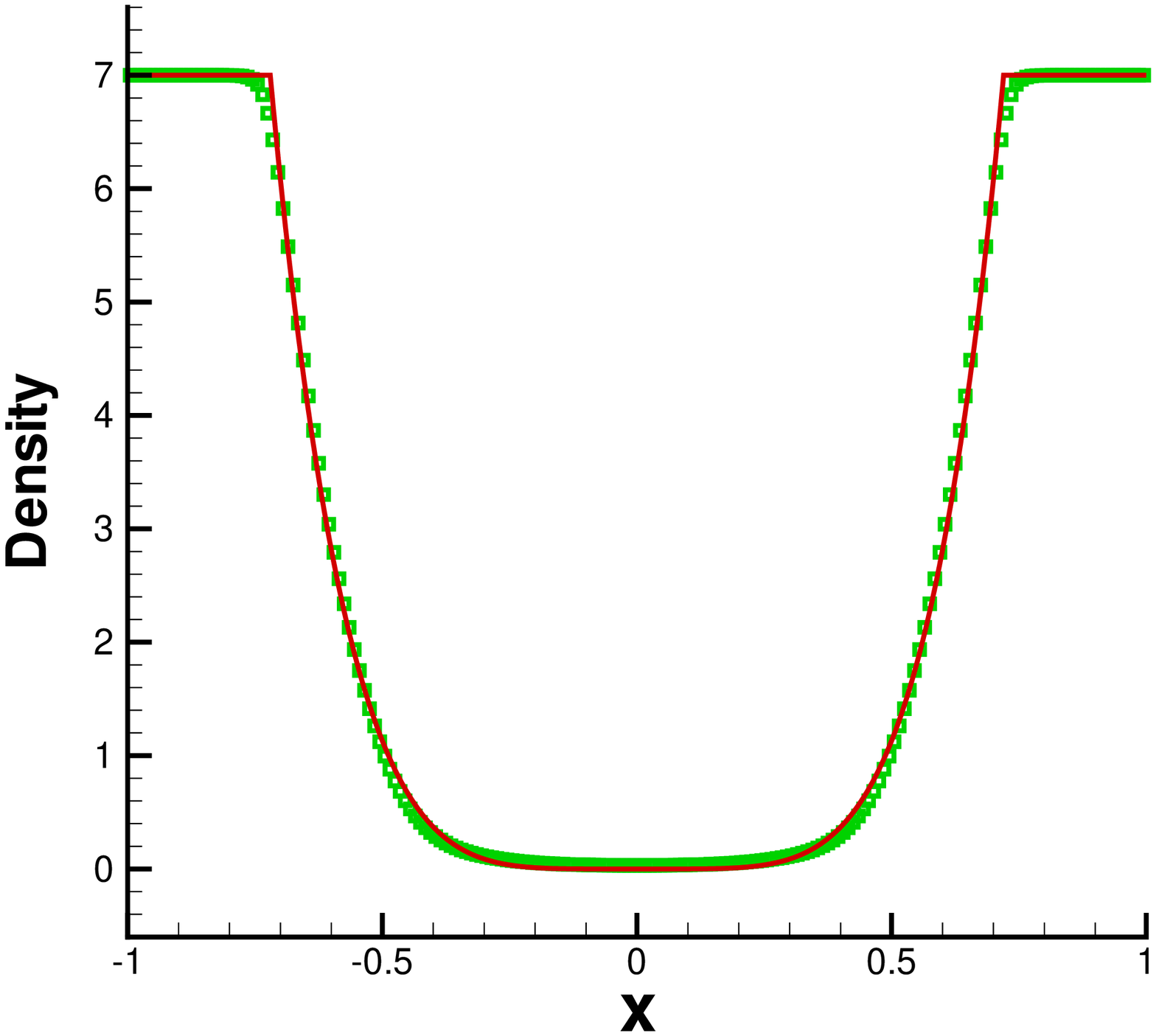},
\includegraphics[totalheight=2.0in]{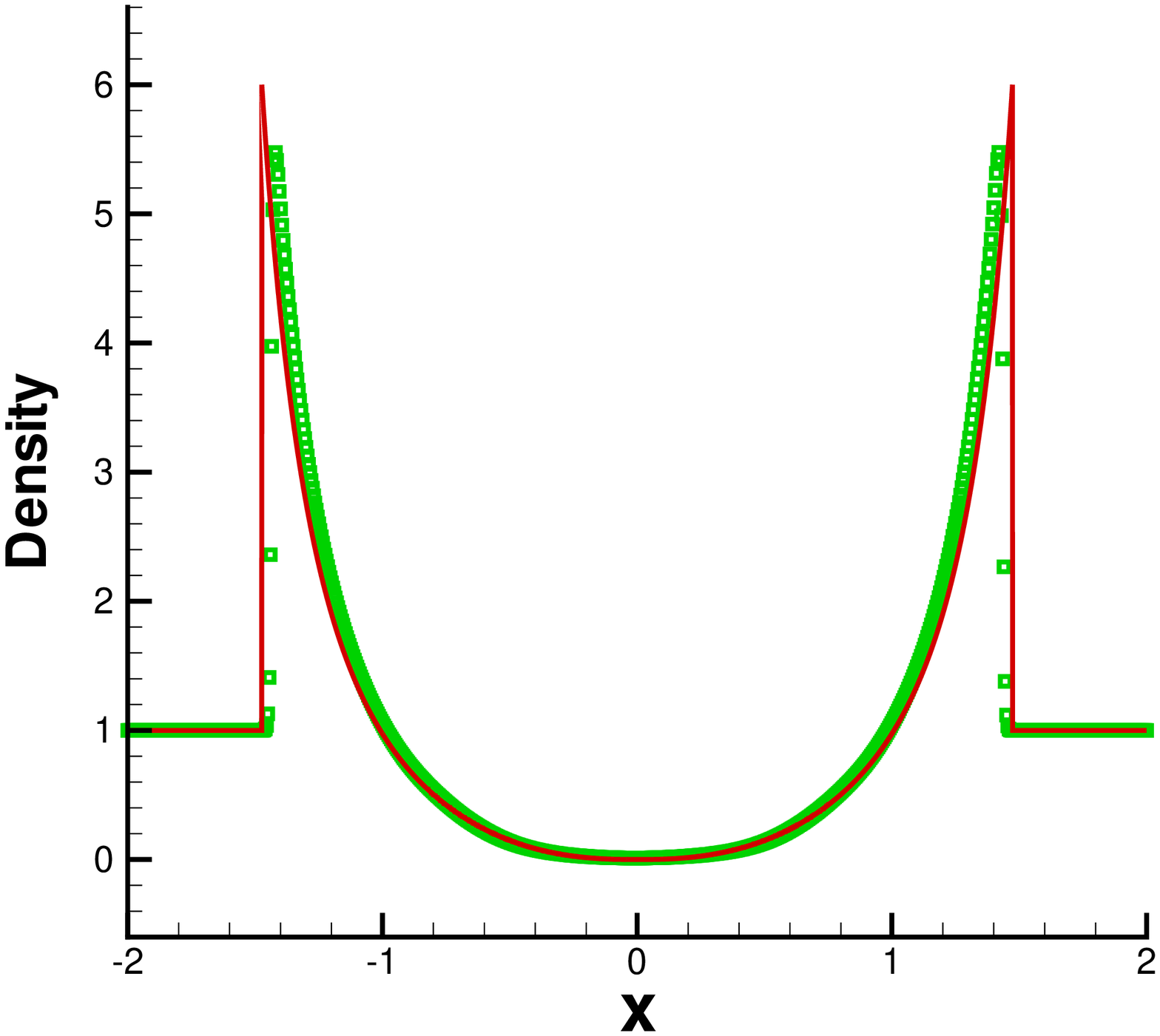}\\
\includegraphics[totalheight=2.0in]{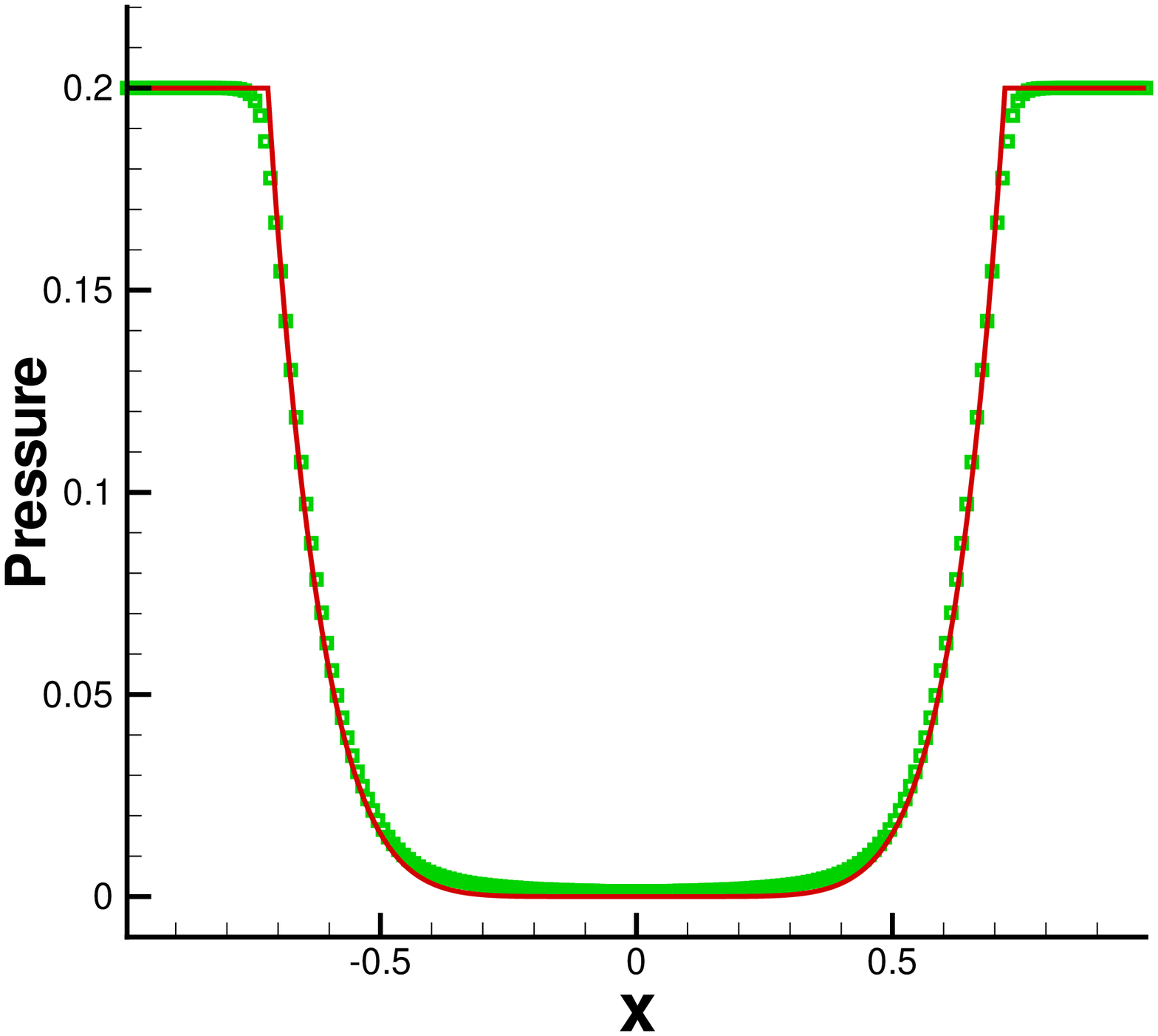},
\includegraphics[totalheight=2.0in]{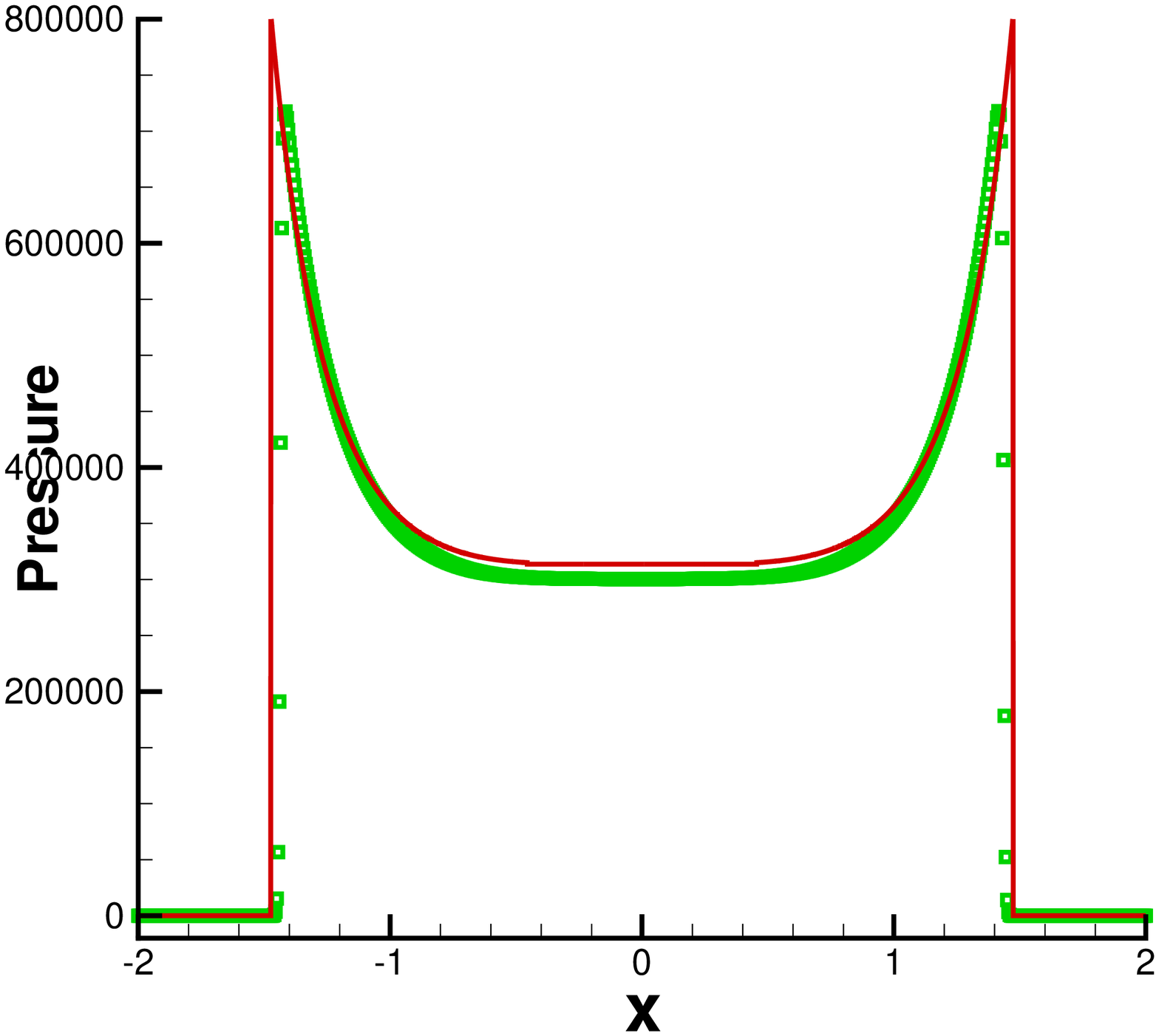}\\
\includegraphics[totalheight=2.0in]{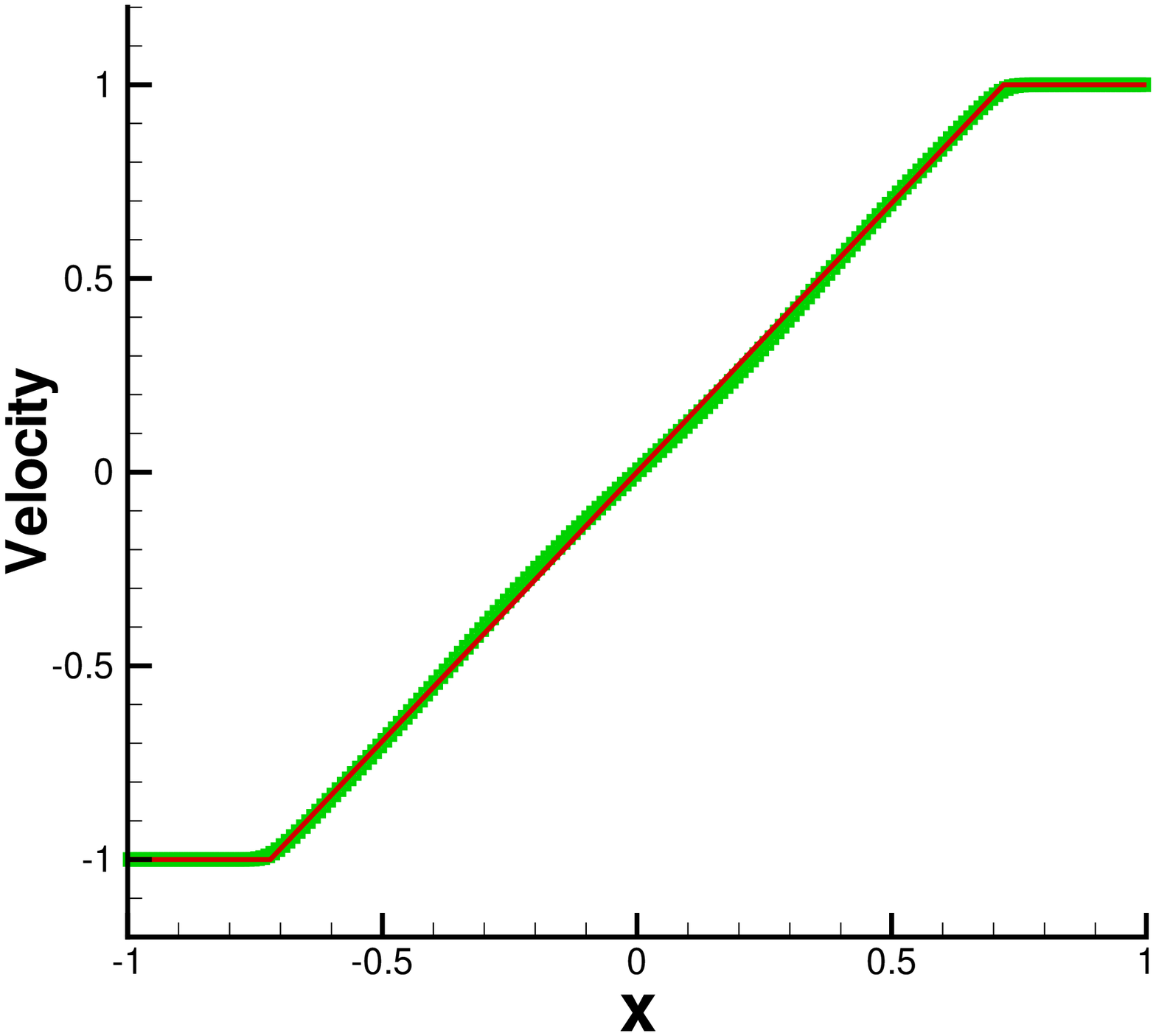},
\includegraphics[totalheight=2.0in]{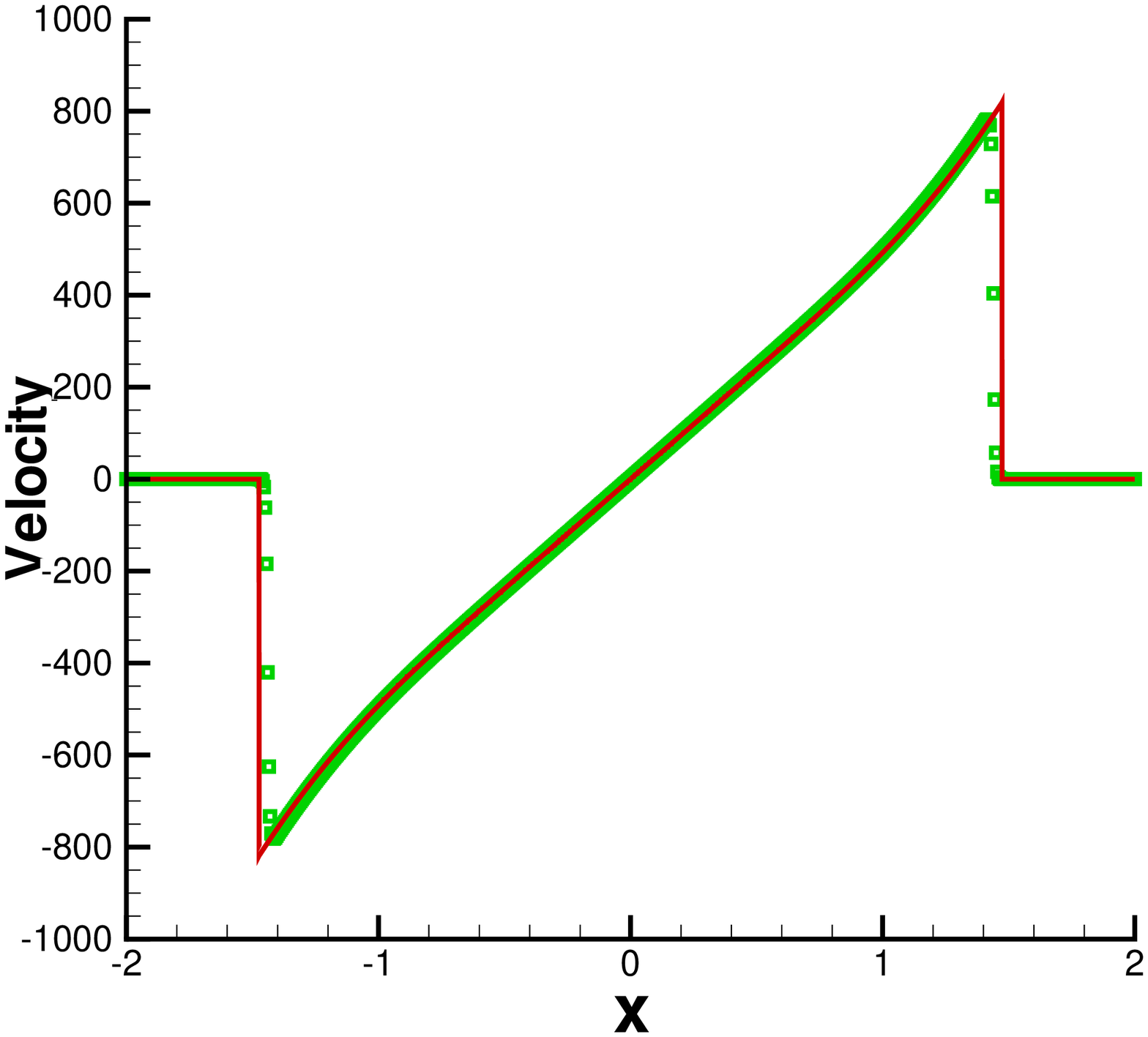}
\caption{Example \ref{ex3}. Left: double rarefaction problem at $T=0.6$. Right: 1D Sedov blast wave at $T=0.001$. $\Delta x=\frac{1}{200}$. The solid lines are the exact solutions. Symbols are the numerical solutions.}
\label{fig52}
\end{figure}

\end{exa}

\begin{exa}{2D low density and low pressure problems.}
\label{ex4}
Now we consider two 2D low density and low pressure problems for the ideal gas.
The first one is the 2D Sedov blast wave. The computational domain is a square of $[0, 1.1]\times[0, 1.1]$.
For the initial condition, similar to the 1D case, the density is $1$, the velocity is $0$, the total energy is $10^{-12}$ everywhere except in the lower left corner is the constant $\frac{0.244816}{\Delta x \Delta y}$.
$\gamma=1.4$. The numerical boundary on the left and bottom edges is reflective. In Fig. \ref{fig53} (left), we show the numerical density at the mesh sizes $\Delta x=\Delta y=\frac{1.1}{160}$ with the PP flux limiters at $T=1$. The numerical solution with cutting along the diagonal matches the exact solution very well in Fig. \ref{fig53} (right).

\begin{figure}
\centering
\includegraphics[totalheight=2.5in,angle=270]{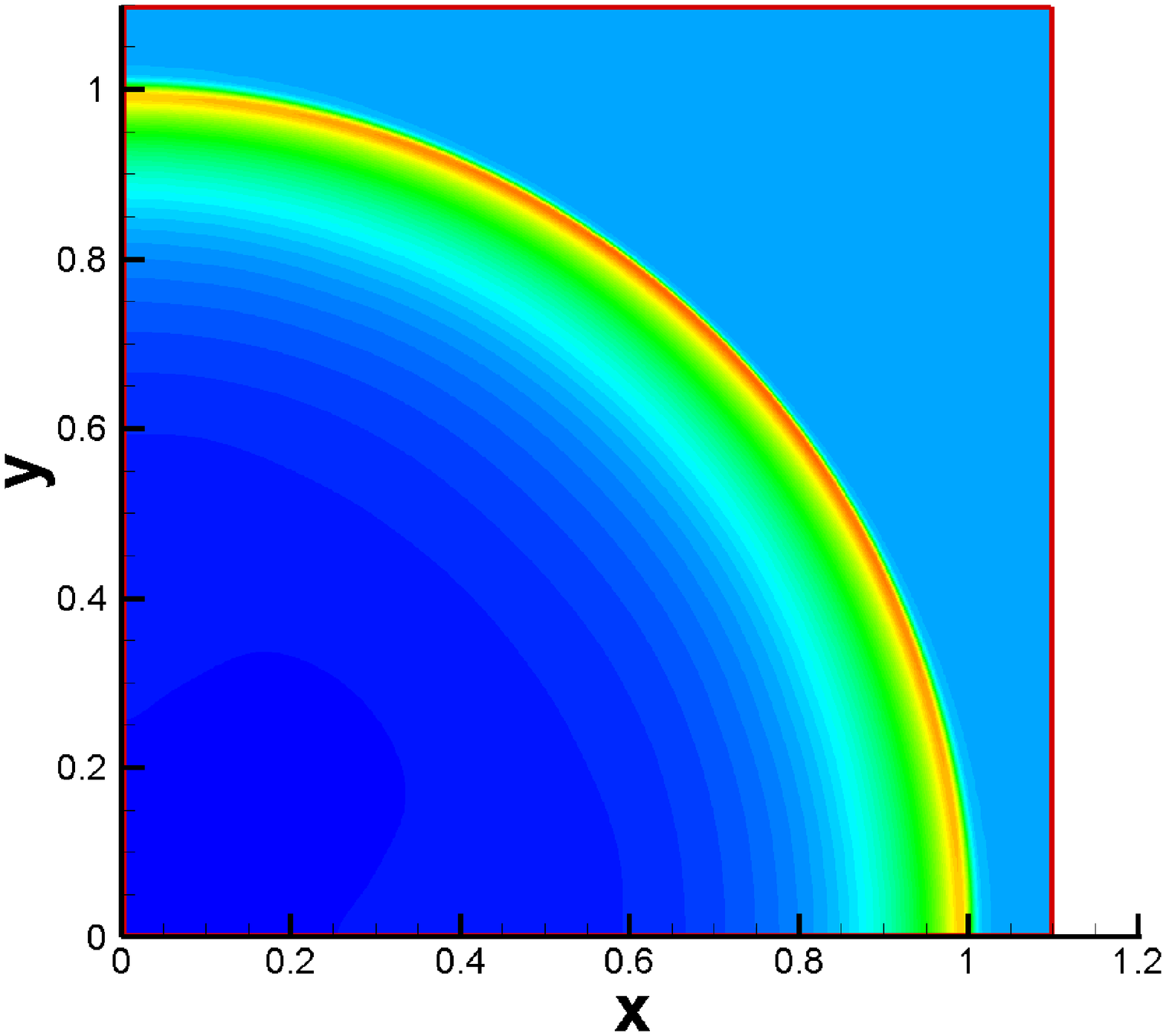}
\includegraphics[totalheight=2.5in,angle=270]{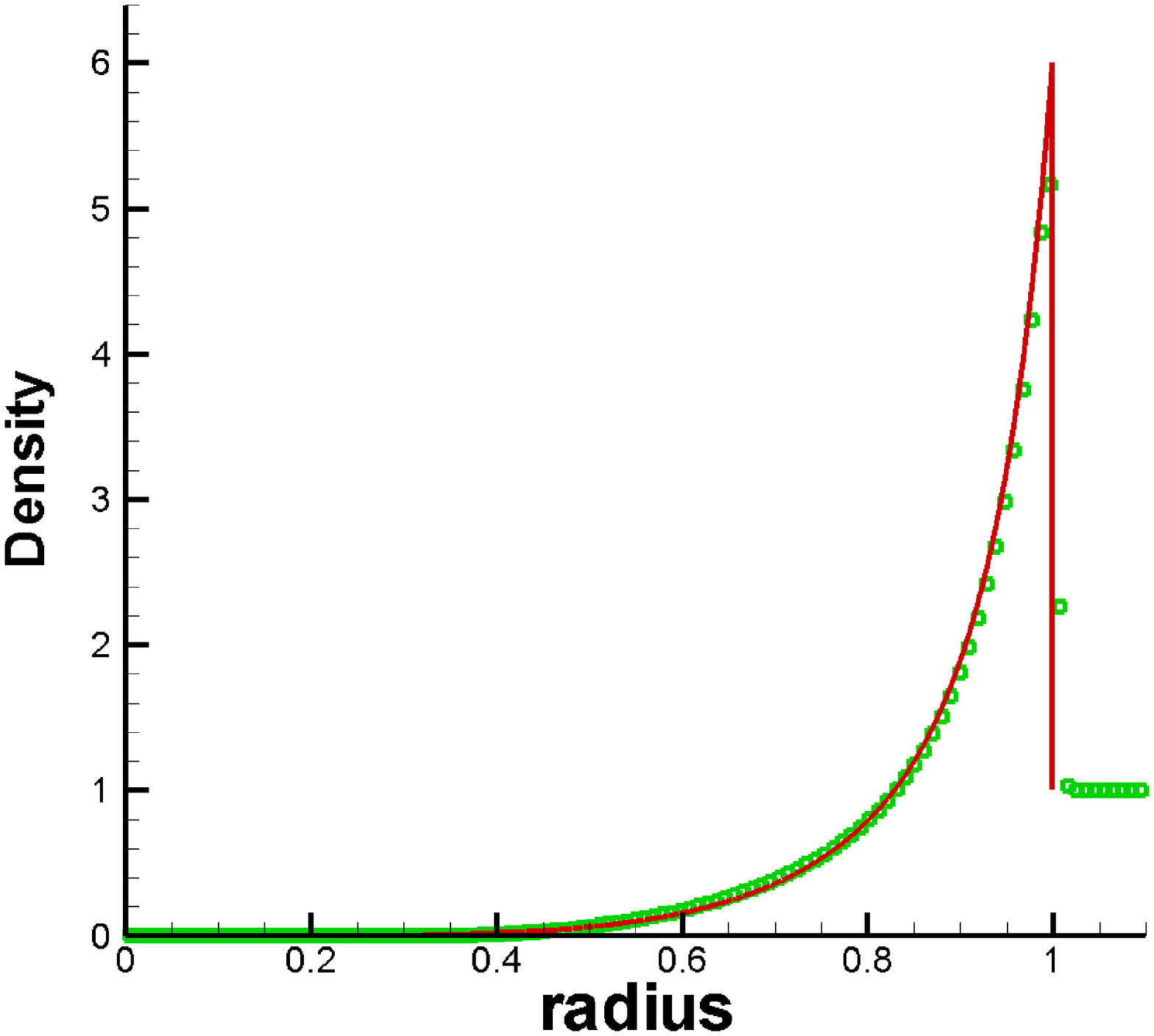}\\
\caption{Example \ref{ex4}. 2D Sedov blast wave. $T=1$. $\Delta x=\Delta y=\frac{1.1}{160}$.
Left: contour of density. Right: cut along diagonal, the solid line is the exact solution, symbols are
the numerical solution.}
\label{fig53}
\end{figure}

The second one is the shock diffraction problem. The computational domain is the union of
$[0,1]\times[6,11]$ and $[1,13]\times[0,11]$. The initial condition is a pure right-moving shock
of $Mach=5.09$, initially located at $x=0.5$ and $6\le y \le 11$, moving into undisturbed air ahead
of the shock. The undisturbed air has a density of $1.4$ and a pressure of $1$. The boundary conditions
are inflow at $x=0$, $6\le y \le 11$, outflow at $x=13$, $0\le y \le 11$, $1\le x\le 13$, $y=0$ and $0\le x \le 13$, $y=11$, and reflective at the walls $0\le x \le 1$, $y=6$ and $x=1$, $0\le y \le 6$. $\gamma=1.4$.
The density and pressure at the mesh sizes $\Delta x=\Delta y=\frac{1}{32}$ with the PP flux limiters at $T=2.3$ are presented in Fig. \ref{fig532}. 

\begin{figure}
\centering
\includegraphics[totalheight=2.5in,angle=270]{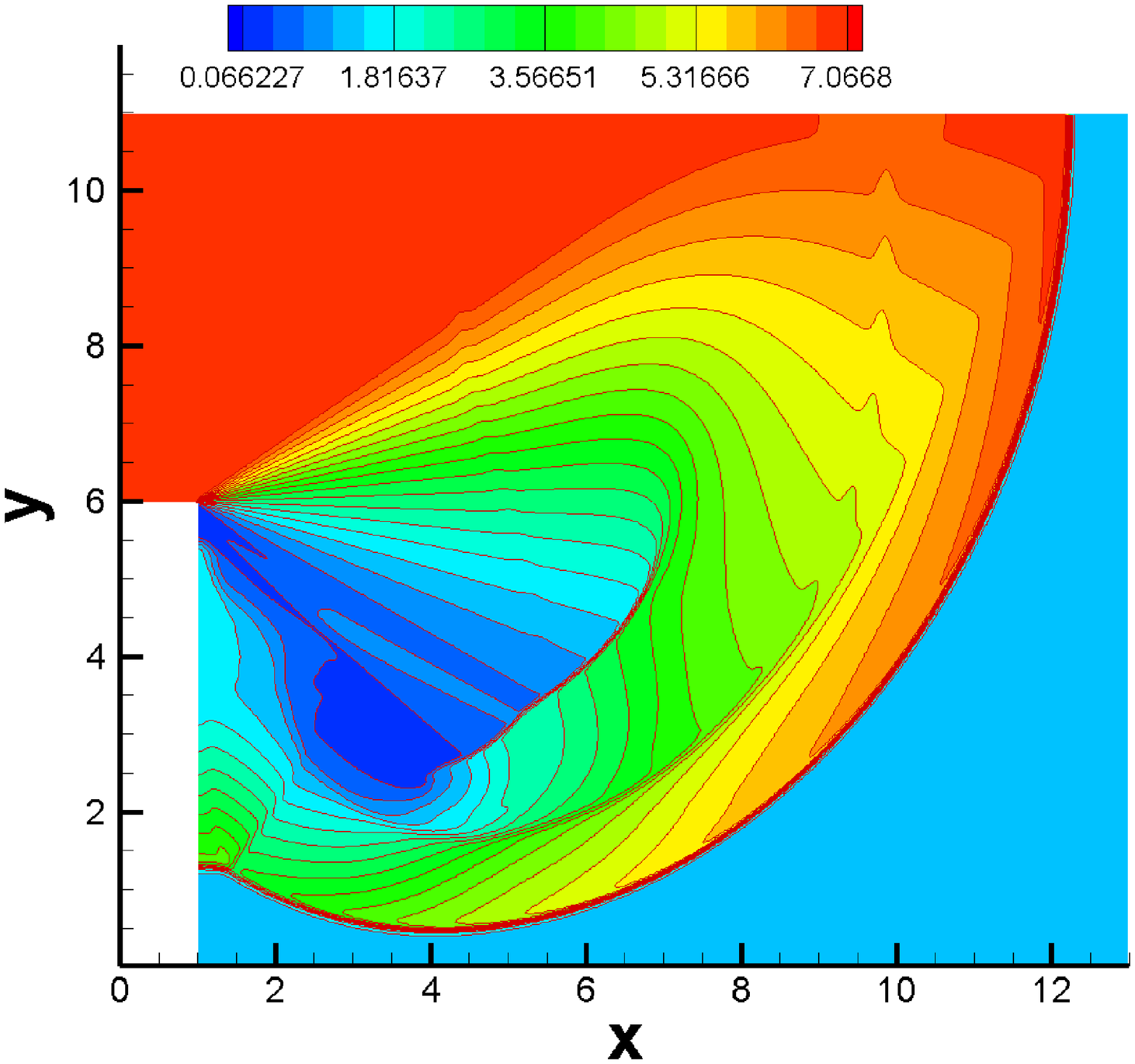},
\includegraphics[totalheight=2.5in,angle=270]{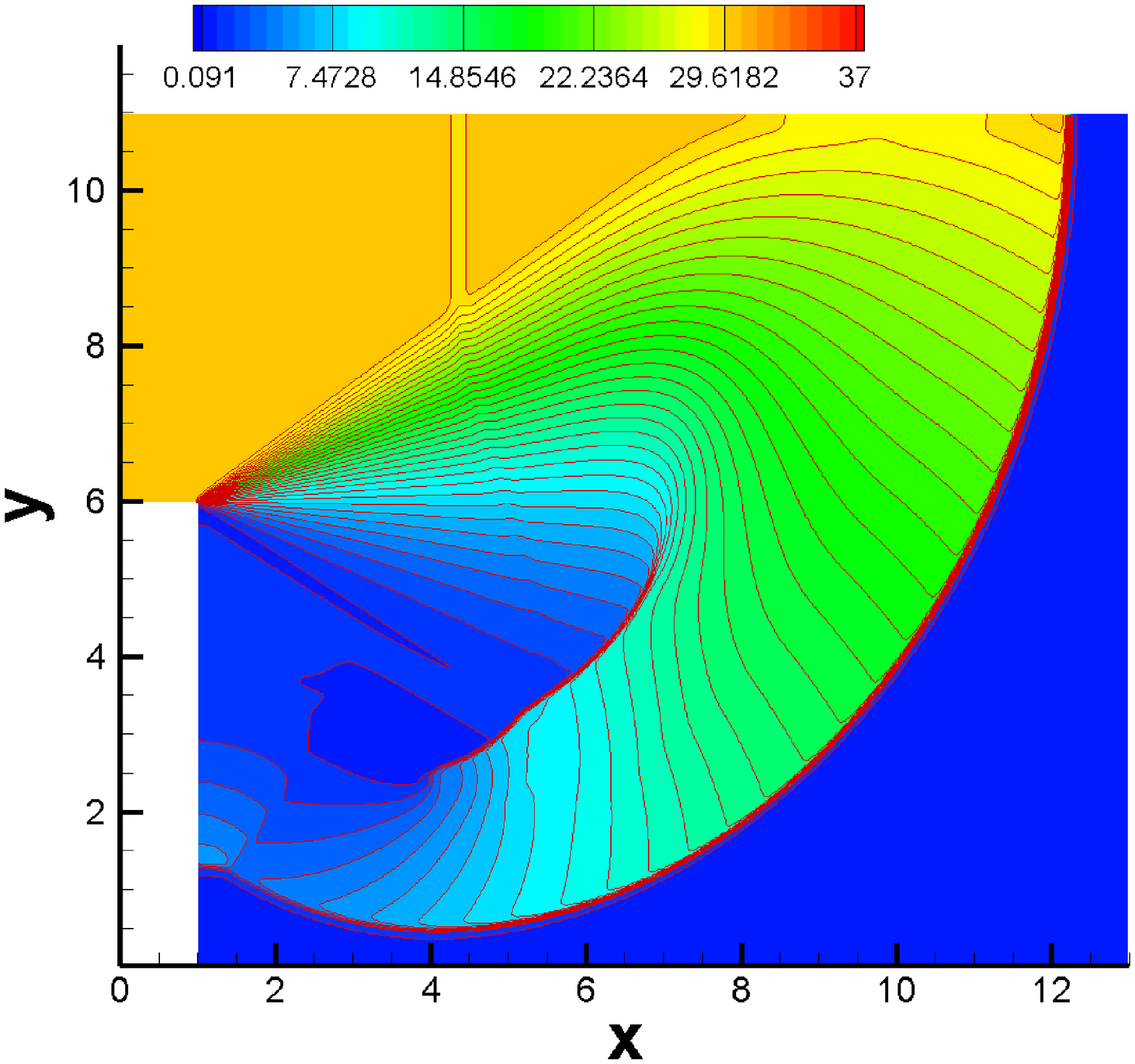}
\caption{Example \ref{ex4}. 2D shock diffraction problem. $T=2.3$. $\Delta x=\Delta y=\frac{1}{32}$.
Left: density, 20 equally spaced contour lines from $\rho=0.066227$ to $\rho=7.0668$. Right: pressure,
40 equally spaced contour lines from $p=0.091$ to $p=37$.}
\label{fig532}
\end{figure}
\end{exa}

\begin{exa}{High Mach number astrophysical jets.}
\label{ex5}
We consider two high Mach number astrophysical jets without the radiative cooling \cite{ha2008positive, zhang2012positivity}.
The first one is a Mach 80 problem. $\gamma=5/3$. The computational domain is $[0,2]\times[-0.5,0.5]$,
which is full of the ambient gas with $(\rho, u, v, p)=(0.5,0,0,0.4127)$ initially. The boundary conditions for the right, top and bottom are outflows. For the left boundary, $(\rho, u, v, p)=(5,30,0,0.4127)$ if $y\in[-0.05, 0.05]$ and $(\rho, u, v, p)=(5,0,0,0.4127)$ otherwise. The numerical density on a mesh of $448\times224$ grid points with the PP flux limiters at $T=0.07$ is shown in Fig. \ref{fig54} (left).
Then a Mach 2000 problem is considered to show the robustness of the scheme with the PP flux limiters. The computational domain is taken as $[0,1]\times[-0.25,0.25]$, initially full of the ambient gas with
$(\rho, u, v, p)=(0.5,0,0,0.4127)$. Similarly, the right, top and bottom boundary are outflows. For the left
boundary,  $(\rho, u, v, p)=(5,800,0,0.4127)$ if $y\in[-0.05, 0.05]$ and $(\rho, u, v, p)=(5,0,0,0.4127)$ otherwise. The numerical density at a mesh of $800\times400$ grid points with the PP flux limiters at $T=0.001$ is shown in Fig. \ref{fig54} (right).

\begin{figure}
\centering
\includegraphics[totalheight=2.5in,angle=270]{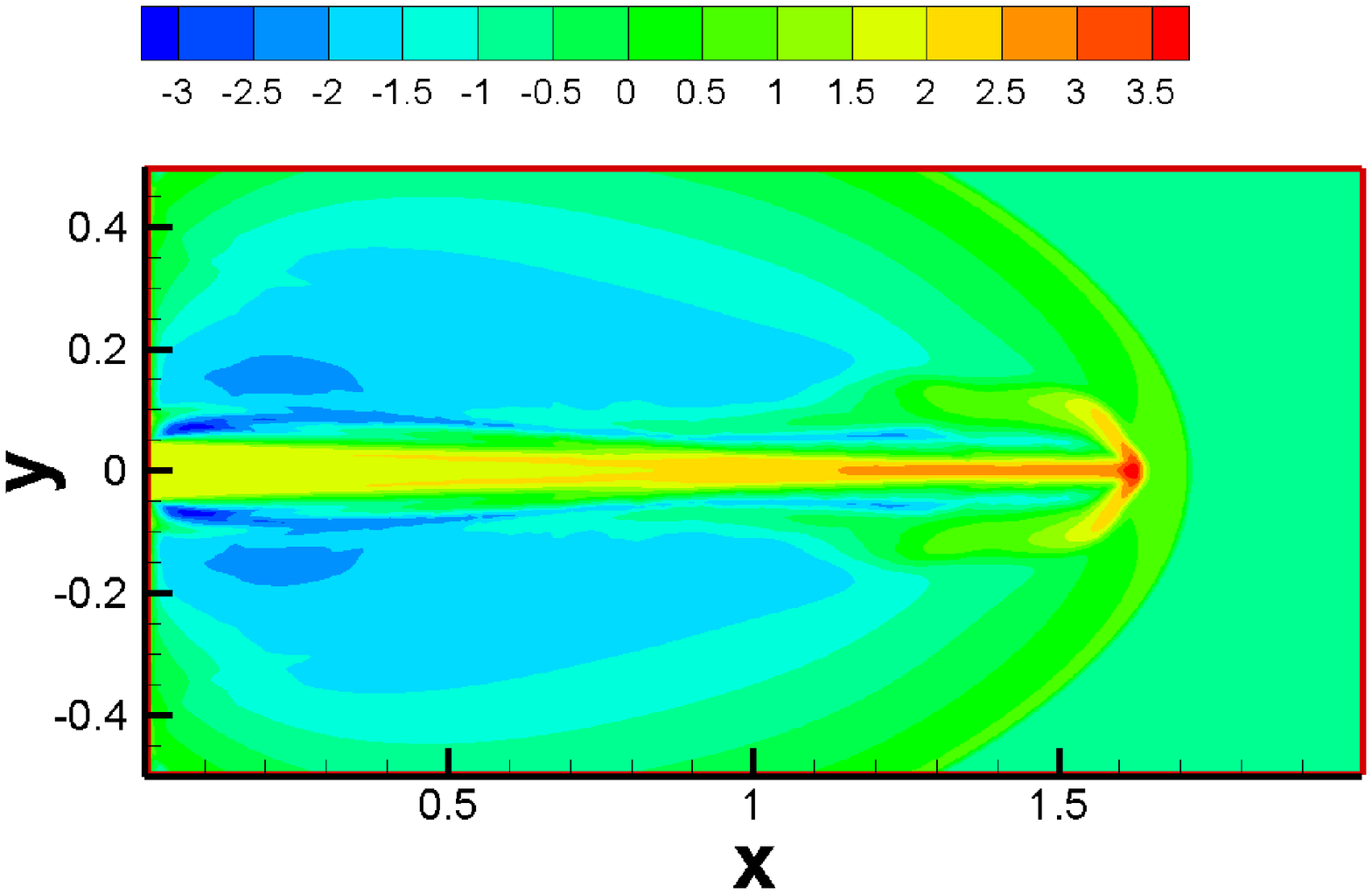}
\includegraphics[totalheight=2.5in,angle=270]{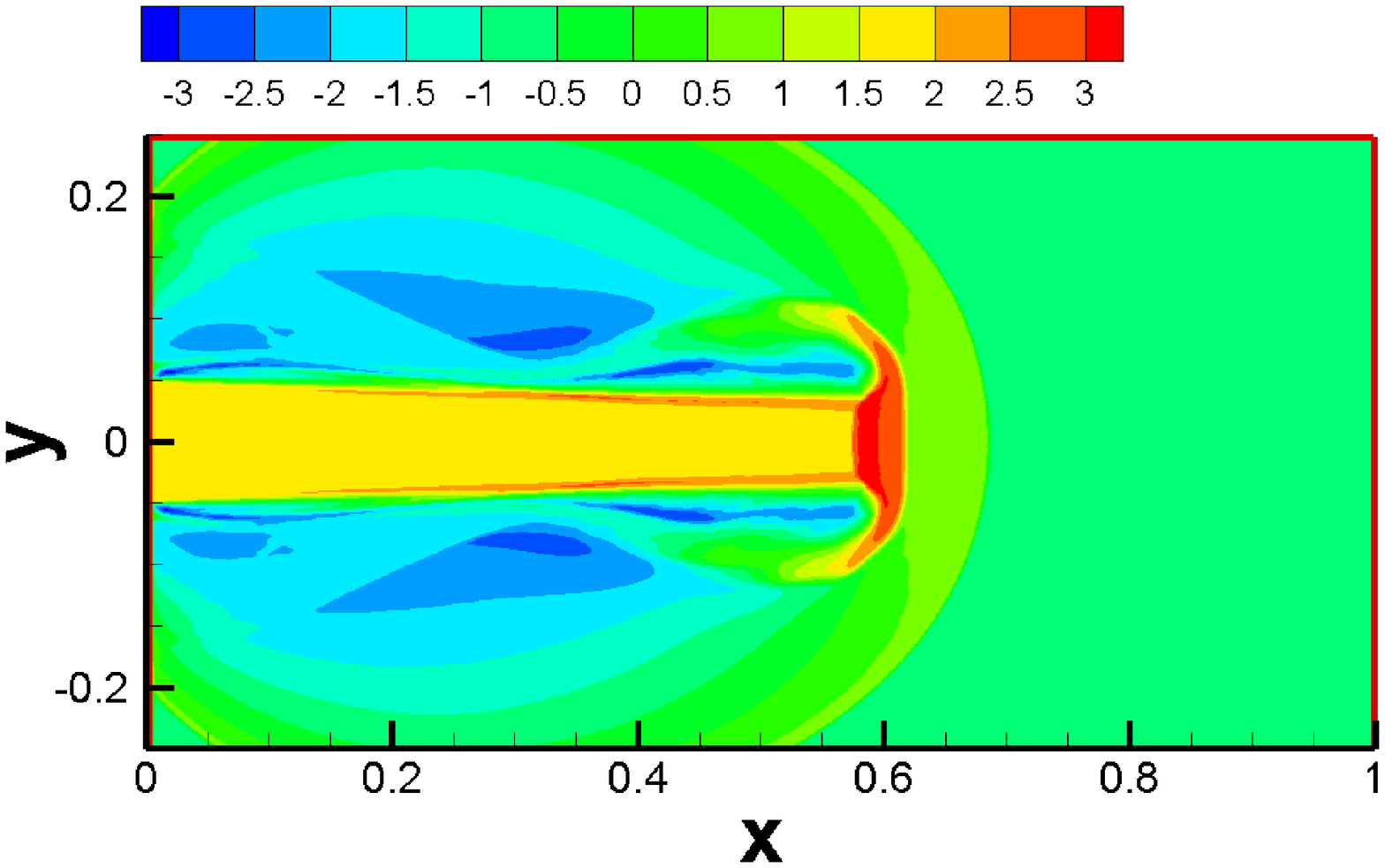}
\caption{Example \ref{ex5}. High Mach number astrophysical jet. $T=2.3$.
Left: density of Mach 80 at $T=0.07$ with mesh $448\times224$.
Right: density of Mach 2000 at $T=0.001$ with mesh $800\times400$.}
\label{fig54}
\end{figure}
\end{exa}

\begin{exa}{The reactive Euler equations.}
\label{ex6}
We consider the following two-dimensional Euler equations with a source term, which are often used to model the detonation waves \cite{wang2011robust, zhang2012positivity}:
\begin{align}
& \mathbf{u}_t+\mathbf{f}(\mathbf{u})_x+\mathbf{g}(\mathbf{u})=\mathbf{s}(\mathbf{u}),
 \quad  t\ge0,\quad (x,y)\in \mathbb{R}^2,
 \label{react1}\\
& \mathbf{u}=
 \begin{pmatrix}
  \rho \\
  m_u \\
  m_v \\
  E \\
  \rho Y
 \end{pmatrix}, \quad
  \mathbf{f}(\mathbf{u})=
 \begin{pmatrix}
  m_u \\
  \rho u^2+p \\
  \rho u v \\
  (E+p)u \\
  \rho u Y
 \end{pmatrix} ,\quad
   \mathbf{g}(\mathbf{u})=
 \begin{pmatrix}
  m_v \\
   \rho u v \\
  \rho v^2+p\\
  (E+p)v \\
  \rho v Y
 \end{pmatrix},\quad
   \mathbf{s}(\mathbf{u})=
 \begin{pmatrix}
  0 \\
 0 \\
  0 \\
  0 \\
  \omega
 \end{pmatrix}
 ,
 \label{react2}
\end{align}
with
\begin{equation*}
 m_u=\rho u,\quad m_v=\rho v,\quad E=\frac{1}{2}\rho u^2+\frac{1}{2}\rho v^2+\frac{p}{\gamma-1}+\rho q Y,
\end{equation*}
where $q$ is the heat release rate of reaction, $\gamma$ is the specific heat ratio and $Y$ is the reactant
mass fraction. The source term is assumed to be in an Arrhenius form
\begin{equation}
 \omega=-\tilde{K}\rho Y\exp(-\tilde{T}/T),
 \label{omega}
\end{equation}
where  $T=\frac{p}{\rho}$ is the temperature, $\tilde{T}$ is the activation temperature and $\tilde{K}$ is a constant.
The eigenvalues of the Jacobian $\mathbf{f}'(\mathbf{u})$ are $u-c, u, u, u, u+c$ and the eigenvalues of the
Jacobian $\mathbf{g}'(\mathbf{u})$ are $v-c, v, v, v, v+c$, where $c=\sqrt{\gamma\frac{p}{\rho}}$.
The computation domain for this problem is the union of $[0,1]\times[2,5]$ and $[1,5]\times[0,5]$. The initial
conditions are, if $x<0.5$,  $(\rho,u,v,E,Y)=(11,6.18,0,970,1)$; otherwise, $(\rho,u,v,E,Y)=(1,0,0,55,1)$. The boundary conditions are reflective except at $x=0$, $(\rho,u,v,E,Y)=(11,6.18,0,970,1)$. Here the parameters are chosen to be $\gamma=1.2$, $q=50$, $\tilde{T}=50$ and $\tilde{K}=2566.4$.

This problem is similar to the shock diffraction problem in Example \ref{ex4}, but this one has a source term. The time step is taken to be
\begin{equation}
 \Delta t=\frac{\text{CFL}}{ \lambda_{max}(\frac{1}{\Delta x}+\frac{1}{\Delta y})+\tilde{K} },
\end{equation}
where $\lambda_{max}=\max\{\||u|+c\|_{\infty}, \||v|+c\|_{\infty}\}$ on all grids, and $\tilde{K}$ comes from the source term (\ref{omega}), such that the first order monotone scheme is PP.
The numerical density and pressure at a mesh of $400\times400$ grid points with the PP flux limiters at $T=0.6$ are shown in Fig. \ref{fig55}, which are comparable to the results in \cite{wang2011robust, zhang2012positivity}.

\begin{figure}
\centering
\includegraphics[totalheight=2.5in,angle=270]{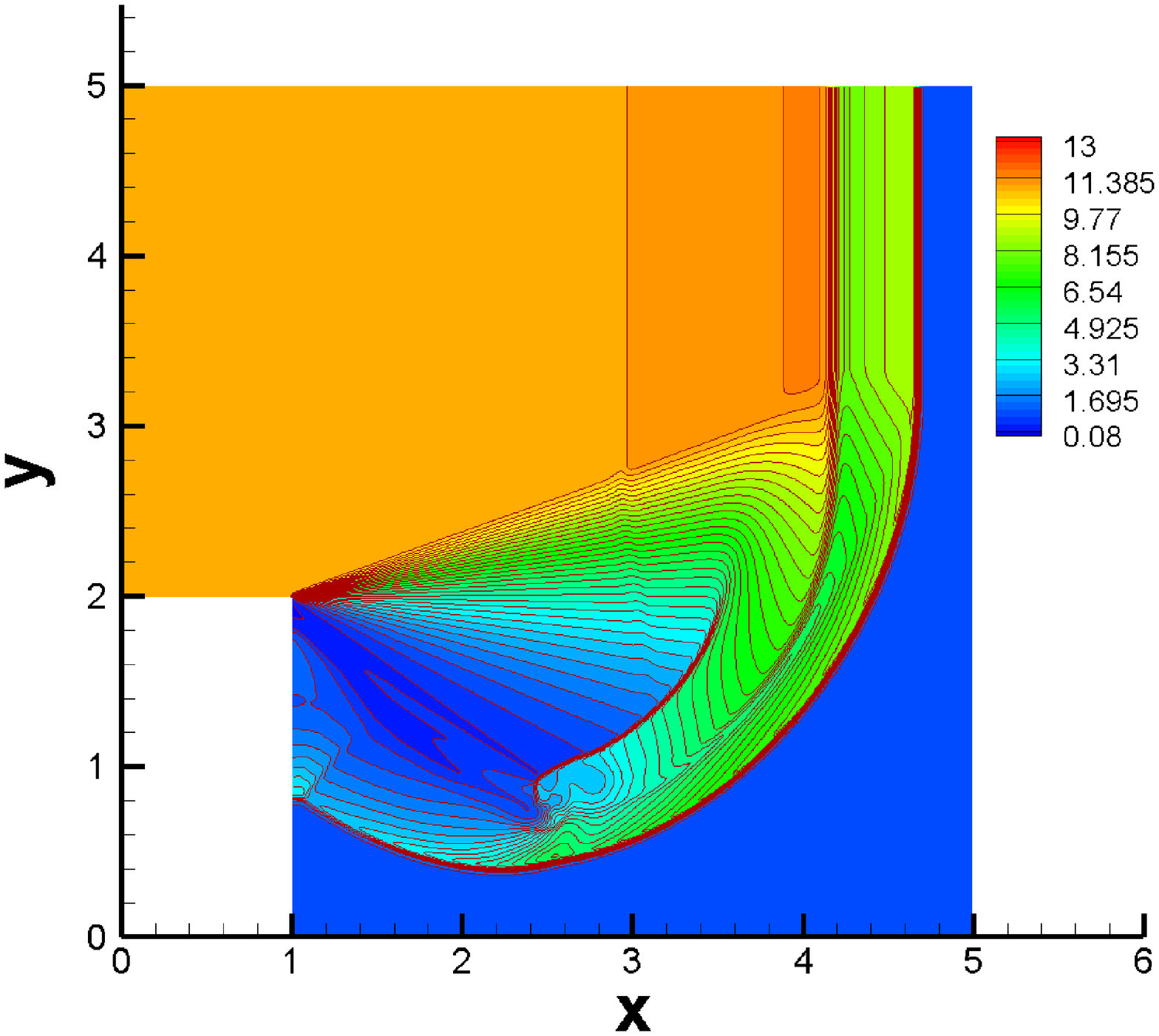}
\includegraphics[totalheight=2.5in,angle=270]{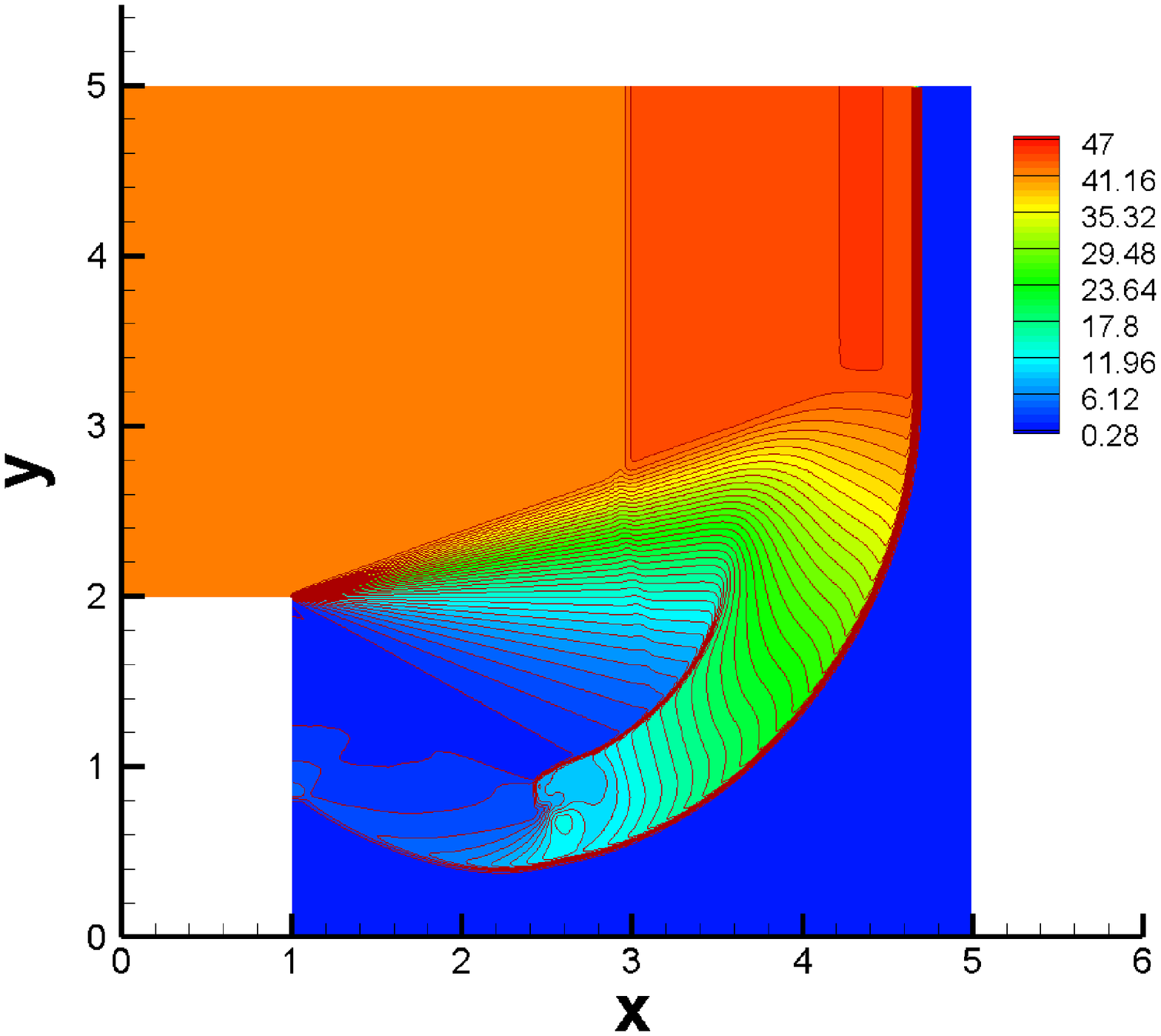}
\caption{Example \ref{ex6}. Detonation diffraction at a $90^\circ$ corner. $T=0.6$.
Mesh $400\times400$. Left: density; Right: pressure.}
\label{fig55}
\end{figure}
\end{exa}

\begin{exa}{General equation of state.}
\label{ex7}
We consider the three species model of the one-dimensional Euler system with a more general equation of state in \cite{wang2009high, zhang2012positivity}.
The model involves three species, $O_2$, $O$ and $N_2$ ($\rho_1=\rho_O$, $\rho_2=\rho_{O_2}$ and $\rho_3=\rho_{N_2}$) with the reaction
\begin{equation}
O_2 + N_2 \rightleftharpoons O + O + N_2.
\end{equation}
The governing equations are
\begin{eqnarray}
\begin{pmatrix}
 \rho_1 \\
 \rho_2 \\
 \rho_3 \\
 \rho u \\
 E
\end{pmatrix}_t
+
\begin{pmatrix}
 \rho_1 u \\
 \rho_2 u \\
 \rho_3 u \\
 \rho u^2 + p \\
 (E+p) u
\end{pmatrix}_x
=
\begin{pmatrix}
 2M_1\omega \\
 -M_2\omega \\
 0 \\
 0 \\
 0
\end{pmatrix}
,
\label{eos}
\end{eqnarray}
and
\begin{equation}
 \rho=\sum_{s=1}^3\rho_s, \quad p=RT\sum_{s=1}^3\frac{\rho_s}{M_s}, \quad E=\sum_{s=1}^3 \rho_s e_s(T)+\rho_1 h_1^0+\frac{1}{2}\rho u^2,
\end{equation}
where the enthalpy $h_1^0$ is a constant, $R$ is the universal gas constant, $M_s$ is the molar mass of species $s$, and the
internal energy $e_s(T)=\frac{3RT}{2M_s}$ and $\frac{5RT}{2M_s}$ for monoatomic and diatomic species respectively. The rate of the chemical reaction is given by
\begin{eqnarray}
 \omega=\left(k_f(T)\frac{\rho_2}{M_2}-k_b(T)\left(\frac{\rho_1}{M_1}\right)^2\right)\sum_{s=1}^3\frac{\rho_s}{M_s}, \quad k_f=C_0 T^{-2}\exp(-E_0/T), \\
 k_b=k_f/\exp(b_1+b_2\log z+ b_3 z+b_4 z^2+b_5 z^3), \quad z=10000/T.
\end{eqnarray}
The parameters and constants are $h_1^0=1.558\times10^7$, $R=8.31447215$, $C_0=2.9\times10^{17} m^3$, $E_0=59750 K$,
and $b_1=2.855$, $b_2=0.988$, $b_3=-6.181$, $b_4=-0.023$, $b_5=-0.001$.

The eigenvalues of the Jacobian are $(u,u,u,u+c,u-c)$ where $c=\sqrt{\gamma\frac{p}{\rho}}$ with $\gamma=1+\frac{p}{T\sum_{s=1}^3\rho_s e_s'(T)}$.
Similar to Example \ref{ex6}, the time step is chosen to be
\begin{equation}
 \Delta t=\frac{\text{CFL } \Delta x}{\lambda_{max}+s_{max} \Delta x},
\end{equation}
here $\lambda_{max}=\max\{\||u|+c\|_{\infty}\}$ on all grids and $s_{max}$ is
\begin{equation}
 s_{max}=\max\left\{\left|\frac{M_2\omega}{\rho_2}\right|,\left|\frac{2M_1\omega}{\rho_1}\right|\right\}.
\end{equation}

A shock tube problem is considered for the reactive flows with high pressure on the left and low pressure on the right initially
in the chemical equilibrium ($\omega=0$). The initial conditions are:
\begin{equation}
 (p_L, T_L)=(1000 N/m^2, 8000K), \quad (p_R, T_R)=(1 N/m^2, 8000K),
\end{equation}
with zero velocity everywhere and the densities satisfying
\begin{equation}
 \frac{\rho_1}{2M_1}+\frac{\rho_2}{M_2}=\frac{21}{79}\frac{\rho_3}{M_3},
\end{equation}
where $M_1=0.016$, $M_2=0.032$ and $M_3=0.028$. The initial densities of $O$, $O_2$ and $N_2$ are
$5.251896311257204\times10^{-5}$, $3.748071704863518\times10^{-5}$ and $2.962489471973072\times10^{-4}$
on the left respectively, and $8.341661837019181\times10^{-8}$,  $9.455418692098664\times10^{-11}$ and $2.748909430004963\times10^{-7}$
on the right respectively.

The numerical solution with the PP flux limiter is computed on a mesh size of $\Delta x=\frac{2}{4000}$ up to $T=0.0001$. $\epsilon_{WENO}=10^{-20}$ is taken as in \cite{zhang2012positivity}. In Fig. \ref{fig56}, the positivity of $\rho_1$, $\rho_2$, $\rho_3$ and $p$ is preserved and converged solutions are observed.

\begin{figure}
\centering
\includegraphics[totalheight=2.0in]{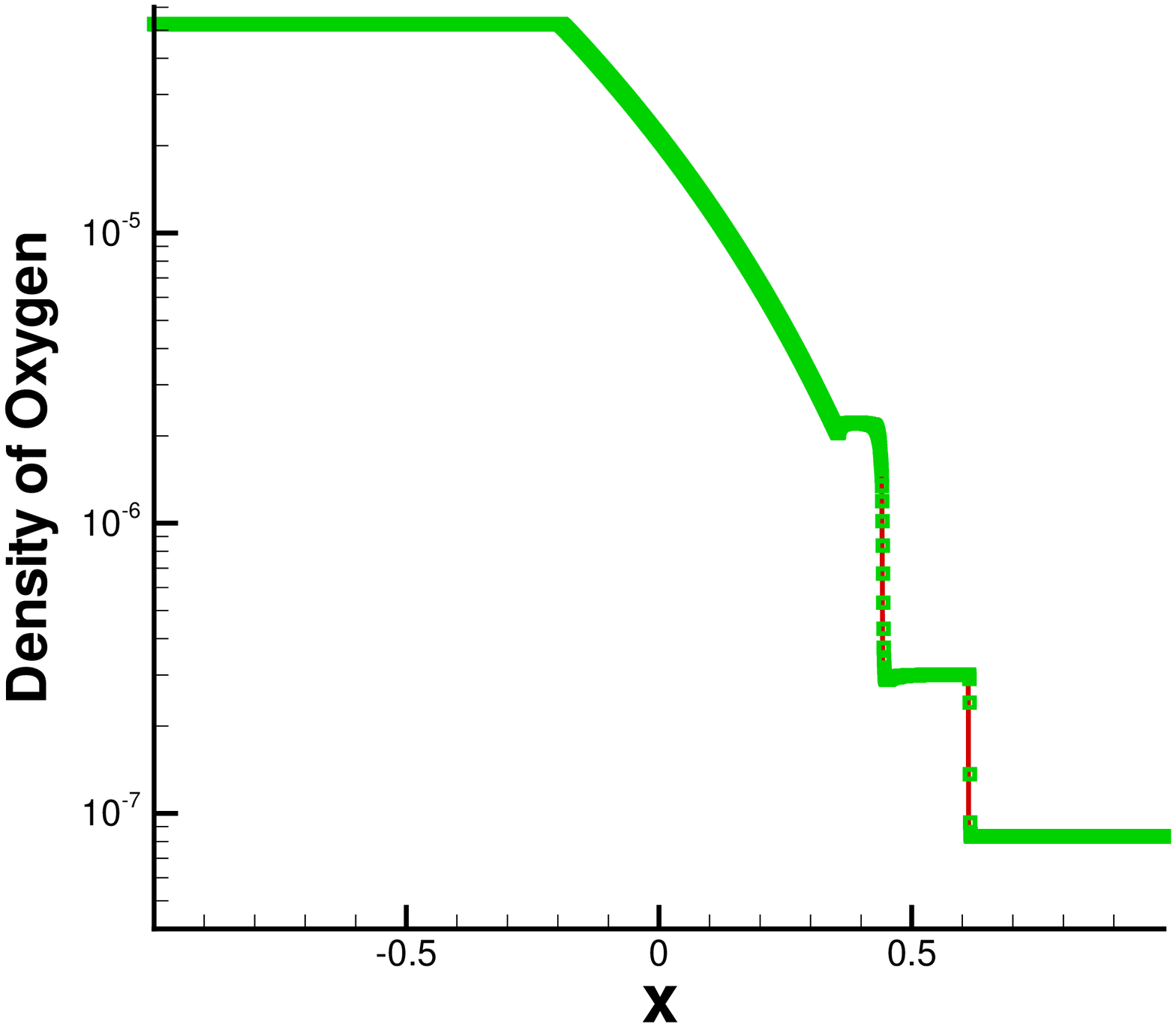},
\includegraphics[totalheight=2.0in]{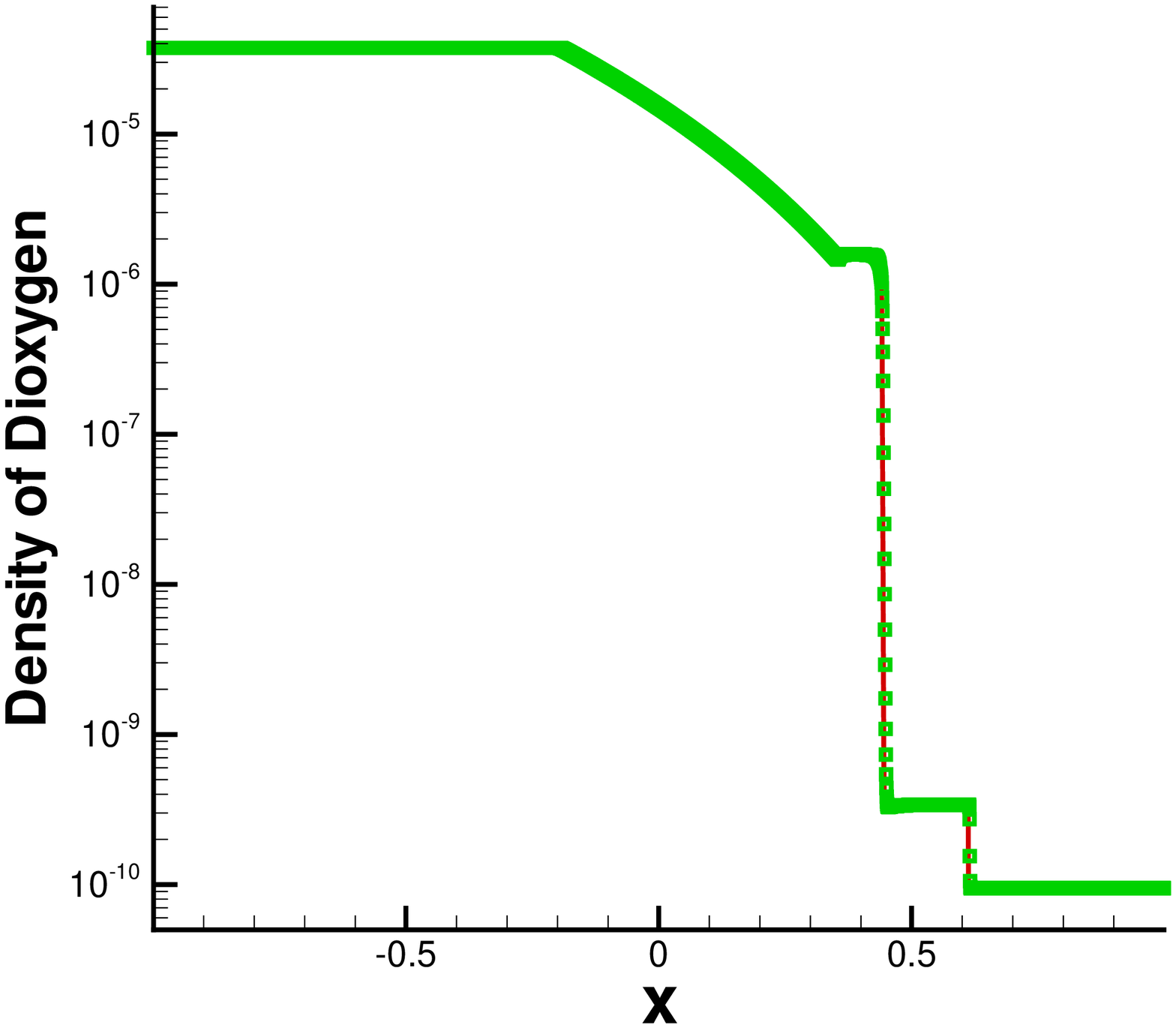}\\
\includegraphics[totalheight=2.0in]{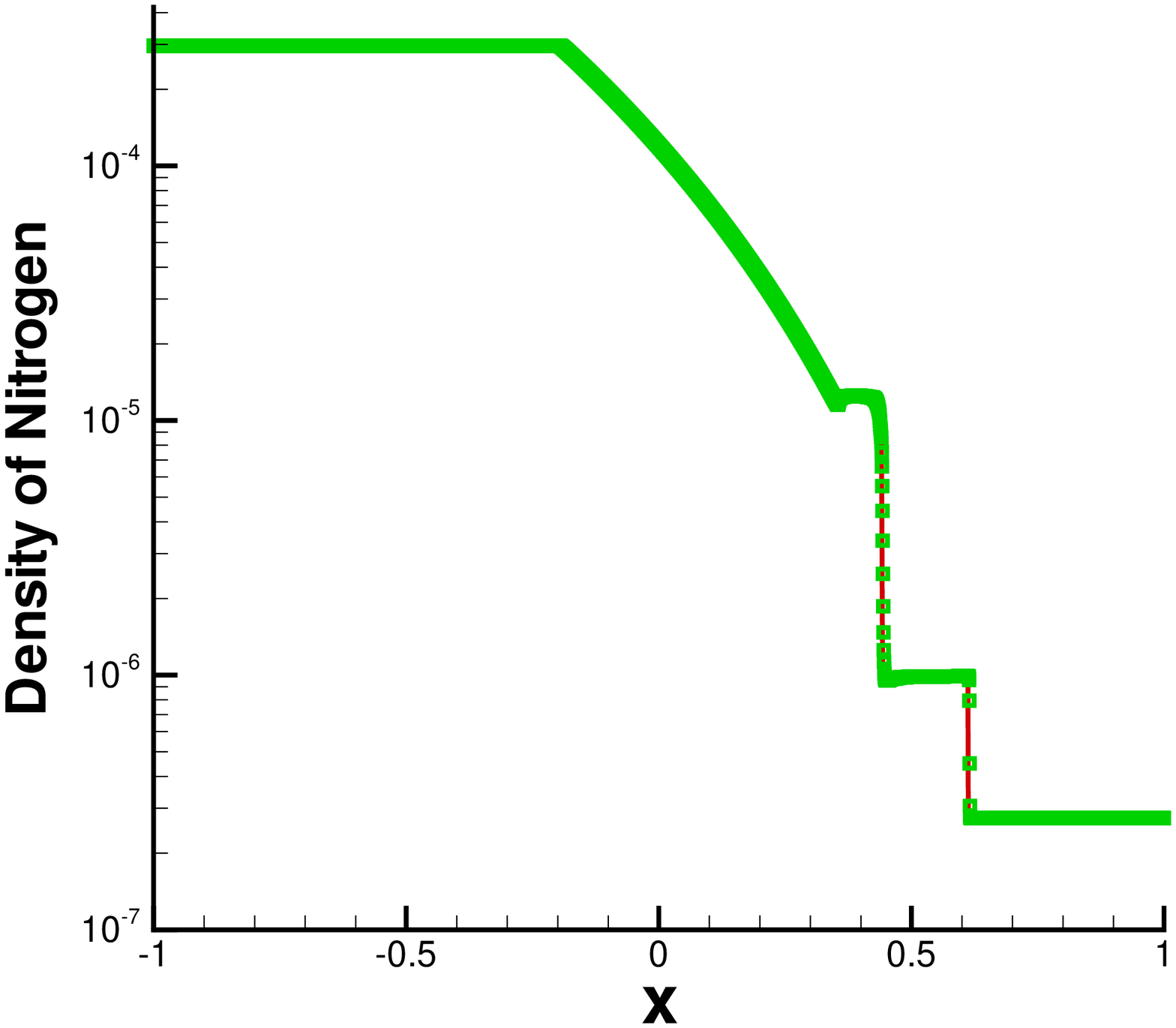},
\includegraphics[totalheight=2.0in]{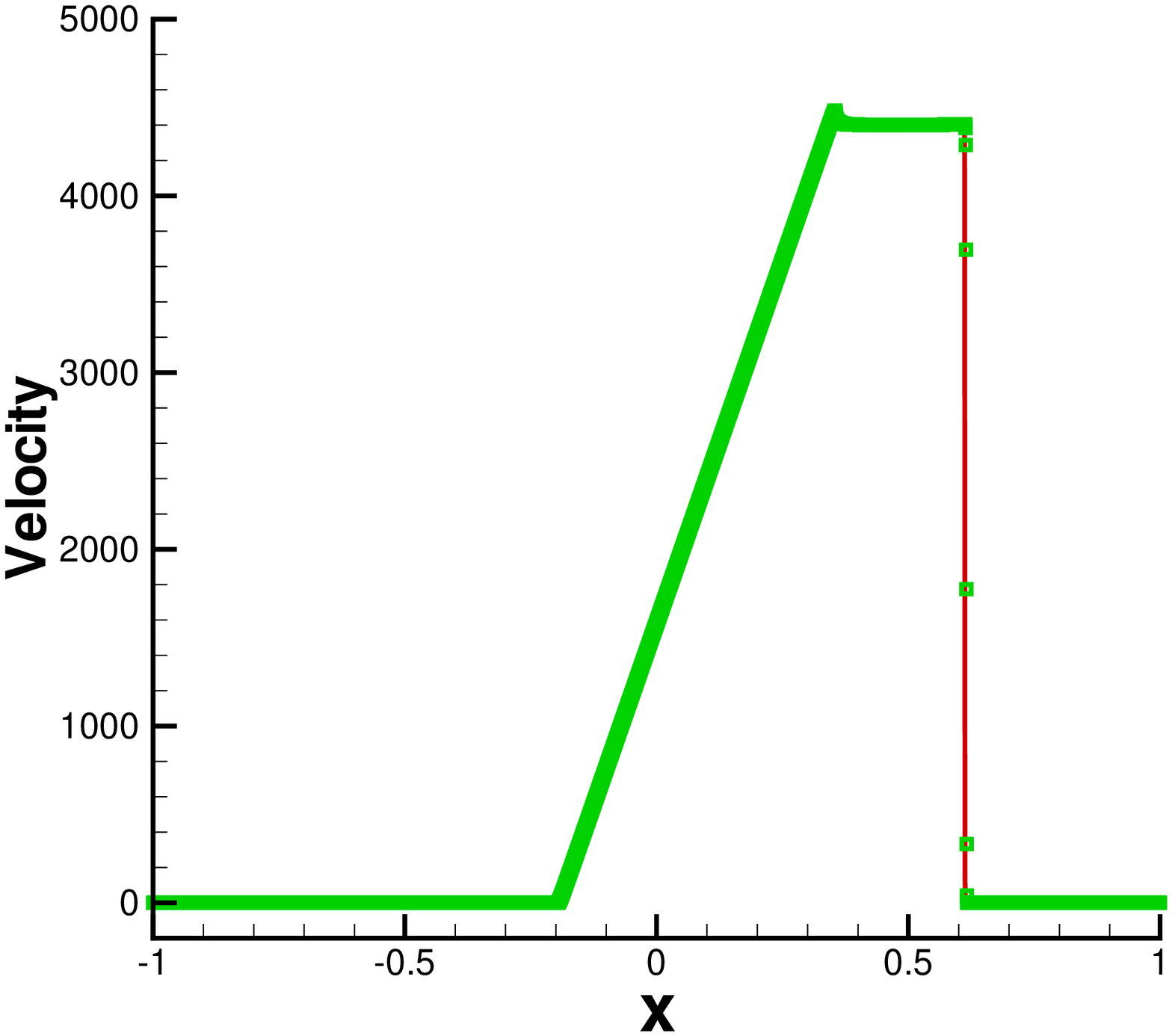}\\
\includegraphics[totalheight=2.0in]{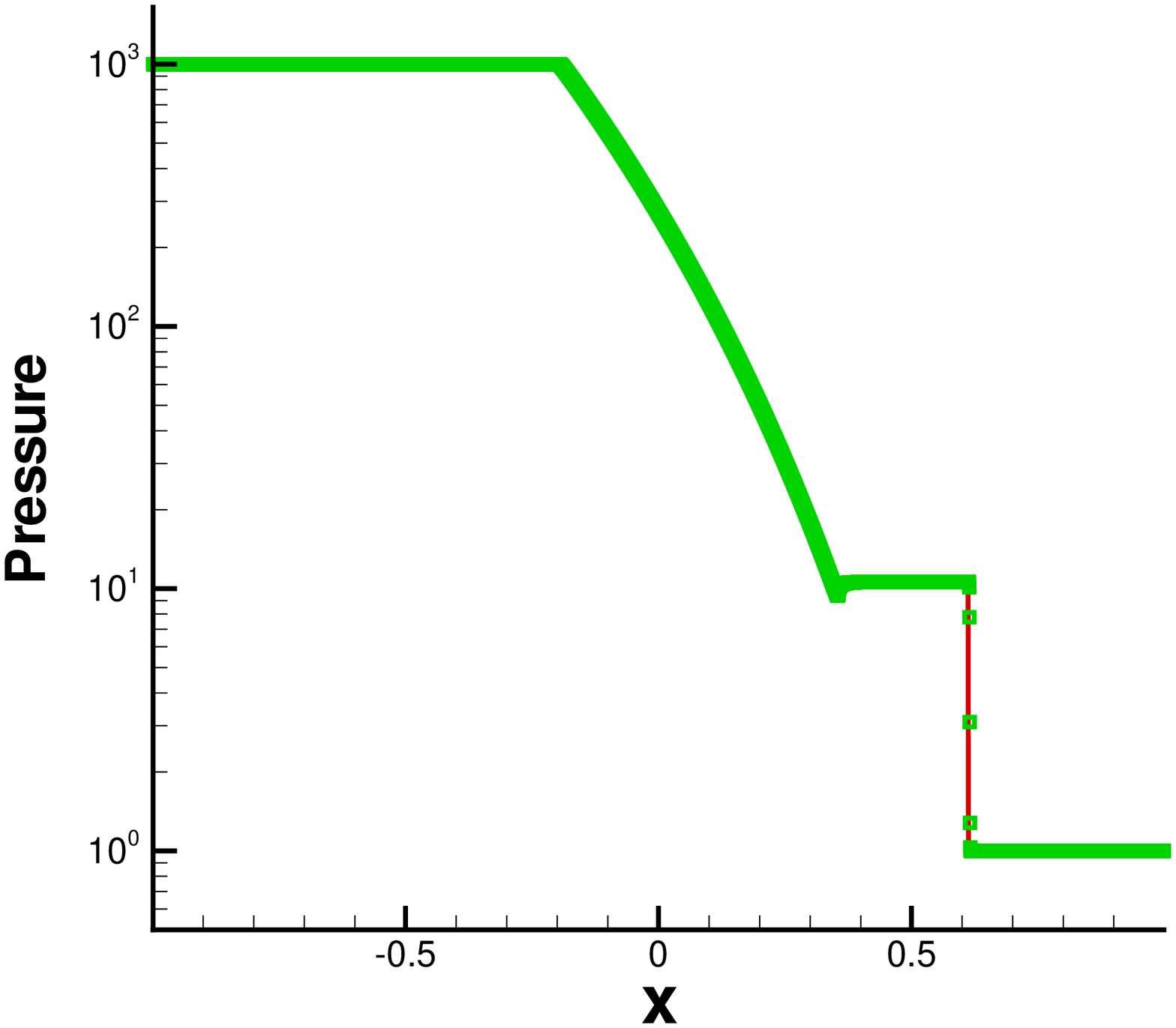}
\caption{Example \ref{ex7}. Three species reaction problem at $T=0.0001$. The
solid lines are the reference solutions at $\Delta x=\frac{2}{8000}$. Symbols are the numerical solutions
at $\Delta x=\frac{2}{4000}$.}
\label{fig56}
\end{figure}
\end{exa}
\section{Conclusion}
\label{sec6}
\setcounter{equation}{0}
\setcounter{figure}{0}
\setcounter{table}{0}
We addressed the potential negative density and pressure problem that emerges when the high order WENO schemes are applied to solve compressible Euler equations in some extreme situations. The approach that we propose is in the conservative high order finite difference WENO approximation framework. We generalized the MPP flux limiting technique for the high order finite difference WENO methods solving scalar conservation law to a class of PP flux limiters for compressible Euler equations. We also developed the parametrized flux limiters for equations with source terms.  Extensive numerical tests show the capability of the proposed approach: without sacrificing accuracy and much of the efficiency, the new schemes produce solutions satisfying the PP property  for scalar problems with a source term, and solutions with positive density and pressure for compressible Euler equations with or without source terms.

\bigskip
\noindent
{\bf Acknowledgement.}
We would like to thank Xiangxiong Zhang from MIT for helpful discussions.

\bibliographystyle{siam}
\bibliography{refer}

\end{document}